
\def\LaTeX{%
  \let\Begin\begin
  \let\End\end
  \let\salta\relax
  \let\finqui\relax
  \let\futuro\relax}

\def\UK{\def\our{our}\let\sz s}
\def\USA{\def\our{or}\let\sz z}

\UK



\LaTeX

\USA


\salta

\documentclass[twoside,12pt]{article}
\setlength{\textheight}{24cm}
\setlength{\textwidth}{16cm}
\setlength{\oddsidemargin}{2mm}
\setlength{\evensidemargin}{2mm}
\setlength{\topmargin}{-15mm}
\parskip2mm


\usepackage[usenames,dvipsnames]{color}
\usepackage{amsmath}
\usepackage{amsthm}
\usepackage{amssymb}
\usepackage[mathcal]{euscript}

\usepackage{cite}
%
%


\definecolor{viola}{rgb}{0.3,0,0.7}
\definecolor{ciclamino}{rgb}{0.5,0,0.5}

\def\gianni #1{{\color{red}#1}}
\def\pier #1{{\color{blue}#1}}
\def\juerg #1{{\color{green}#1}}

\def\pier #1{#1}
\def\juerg #1{#1}
\def\gianni #1{#1}




\bibliographystyle{plain}


%

\finqui

\def\Beq{\Begin{equation}}
\def\Eeq{\End{equation}}
\def\Bsist{\Begin{eqnarray}}
\def\Esist{\End{eqnarray}}

\def\Bthm{\Begin{theorem}}
\def\Ethm{\End{theorem}}
\def\Blem{\Begin{lemma}}
\def\Elem{\End{lemma}}
\def\Bprop{\Begin{proposition}}
\def\Eprop{\End{proposition}}
\def\Bcor{\Begin{corollary}}
\def\Ecor{\End{corollary}}
\def\Brem{\Begin{remark}\rm}
\def\Erem{\End{remark}}

\def\Bdim{\Begin{proof}}
\def\Edim{\End{proof}}
\def\Bcenter{\Begin{center}}
\def\Ecenter{\End{center}}
\let\non\nonumber




\def\step #1 \par{\medskip\noindent{\bf #1.}\quad}


\def\aand{\quad\hbox{and}\quad}

\def\lhs{left-hand side}
\def\rhs{right-hand side}



\def\multibold #1{\def\arg{#1}%
  \ifx\arg\pto \let\next\relax
  \else
  \def\next{\expandafter
    \def\csname #1#1#1\endcsname{{\bf #1}}%
    \multibold}%
  \fi \next}

\def\pto{.}

\def\multical #1{\def\arg{#1}%
  \ifx\arg\pto \let\next\relax
  \else
  \def\next{\expandafter
    \def\csname cal#1\endcsname{{\cal #1}}%
    \multical}%
  \fi \next}


\def\multimathop #1 {\def\arg{#1}%
  \ifx\arg\pto \let\next\relax
  \else
  \def\next{\expandafter
    \def\csname #1\endcsname{\mathop{\rm #1}\nolimits}%
    \multimathop}%
  \fi \next}

\multibold
qwertyuiopasdfghjklzxcvbnmQWERTYUIOPASDFGHJKLZXCVBNM.

\multical
QWERTYUIOPASDFGHJKLZXCVBNM.

\multimathop
diag dist div dom mean meas sign supp .

\def\infess{\mathop{\rm inf\,ess}}
\def\supess{\mathop{\rm sup\,ess}}


\def\accorpa #1#2{\eqref{#1}--\eqref{#2}}
\def\Accorpa #1#2 #3 {\gdef #1{\eqref{#2}--\eqref{#3}}%
  \wlog{}\wlog{\string #1 -> #2 - #3}\wlog{}}


\def\somma #1#2#3{\sum_{#1=#2}^{#3}}

\def\<#1>{\mathopen\langle #1\mathclose\rangle}
\def\norma #1{\mathopen \| #1\mathclose \|}

\def\[#1]{\mathopen\langle\!\langle #1\mathclose\rangle\!\rangle}

\def\ioT {\int_0^T}

\def\intQ{\int_Q}
\def\iO{\int_\Omega}

\def\dt{\partial_t}
\def\dn{\partial_{\bf n}}
\def\indi{I_{[-1,1]}}
\def\dI{\partial I_{[-1,1]}} 

\def\cpto{\,\cdot\,}

\def\checkmmode #1{\relax\ifmmode\hbox{#1}\else{#1}\fi}


\def\erre{{\mathbb{R}}}

\def\enne{{\mathbb{N}}}




\def\genspazio #1#2#3#4#5{#1^{#2}(#5,#4;#3)}
\def\spazio #1#2#3{\genspazio {#1}{#2}{#3}T0}

\def\L {\spazio L}


\def\Lx #1{L^{#1}(\Omega)}
\def\Hx #1{H^{#1}(\Omega)}

\def\Luno{\Lx 1}
\def\Ldue{\Lx 2}

\def\Huno{\Hx 1}
\def\Hdue{\Hx 2}
\def\Hunoz{{H^1_0(\Omega)}}



\let\theta\vartheta

\let\phi\varphi

\let\hat\widehat

\let\TeXchi\chi                         
\newbox\chibox
\setbox0 \hbox{\mathsurround0pt $\TeXchi$}
\setbox\chibox \hbox{\raise\dp0 \box 0 }
\def\chi{\copy\chibox}



\def\Az #1{A_0^{#1}}
\def\Vz #1{V_0^{#1}}
\def\VA #1{V_A^{#1}}
\def\VB #1{V_B^{#1}}

\def\yal{y^{\alpha}}
\def\mual{\mu^{\alpha}}
\def\yn{y^{\alpha_n}}
\def\mn{\mu^{\alpha_n}}
\def\un{u^{\alpha_n}}

\def\bu{\bar u}
\def\by{\bar y}
\def\bm{\bar \mu}
\def\bual{\bar u^{\alpha}}
\def\byal{\bar y^{\alpha}}
\def\bmal{\bar \mu^{\alpha}}
\def\qal{q^\alpha}
\def\pal{p^\alpha}

\def\CP0{(${\mathcal{CP}}_0$)}
\def\CPal{(${\mathcal{CP}}_\alpha$)}

\def\S{\mathcal{S}}
\def\uad{{\mathcal{U}}_{\rm ad}}

\Begin{document}


%
\title{Deep quench approximation and optimal control \\of general Cahn--Hilliard systems
with fractional \\operators and double obstacle potentials}

\author{}
\date{}
\maketitle
\Bcenter
\vskip-1cm
{\large\sc Pierluigi Colli$^{(1)}$}\\
{\normalsize e-mail: {\tt pierluigi.colli@unipv.it}}\\[.25cm]
{\large\sc Gianni Gilardi$^{(1)}$}\\
{\normalsize e-mail: {\tt gianni.gilardi@unipv.it}}\\[.25cm]
{\large\sc J\"urgen Sprekels$^{(2)}$}\\
{\normalsize e-mail: {\tt sprekels@wias-berlin.de}}\\[.45cm]
$^{(1)}$
{\small Dipartimento di Matematica ``F. Casorati'', Universit\`a di Pavia}\\
{\small and Research Associate at the IMATI -- C.N.R. Pavia}\\
{\small via Ferrata 5, 27100 Pavia, Italy}\\[.2cm]
$^{(2)}$
{\small Department of Mathematics}\\
{\small Humboldt-Universit\"at zu Berlin}\\
{\small Unter den Linden 6, 10099 Berlin, Germany}\\[2mm]
{\small and}\\[2mm]
{\small Weierstrass Institute for Applied Analysis and Stochastics}\\
{\small Mohrenstrasse 39, 10117 Berlin, Germany}
\Ecenter
\Begin{abstract}\noindent
In the recent paper \gianni{{\sl ``Well-posedness and regularity for a generalized fractional Cahn--Hilliard system''},
the same authors derived general well-posedness and regularity results} 
for a rather general system of evolutionary operator equations having the structure of a  Cahn--Hilliard system. 
The operators  appearing in the system equations were fractional versions in the spectral sense of general 
linear operators $A$ and $B$ having
compact resolvents and are densely defined,  unbounded, selfadjoint, and monotone
in a Hilbert space of functions defined in a smooth domain. 
The associated double-well potentials driving the phase separation process modeled by the Cahn--Hilliard system 
could be of a very general type that includes standard physically meaningful cases such as 
polynomial, logarithmic, and double obstacle nonlinearities. 
In the \pier{subsequent} paper \gianni{{\sl ``Optimal distributed control of a generalized 
fractional Cahn--Hilliard system''} \juerg{(Appl. Math. Optim. (2018), 
https://doi.org/10.1007/s00245-018-9540-7)} by the same authors}, 
an analysis of distributed optimal control problems was performed for such
evolutionary systems, where only the differentiable case of certain polynomial
and logarithmic double-well potentials could be admitted. Results concerning existence of optimizers
and first-order necessary optimality conditions were derived, where more restrictive 
conditions on the operators $A$ and $B$ had to be assumed in order to be able to show differentiability
properties for the associated control-to-state operator. 
In the present paper, we complement \gianni{these results} by studying a distributed control problem for
such evolutionary systems in the case of nondifferentiable nonlinearities of double obstacle type. 
For such nonlinearities, it is well known that the standard constraint qualifications cannot be
applied to construct appropriate Lagrange multipliers. 
To overcome this difficulty, we follow here the so-called ``deep quench'' method.  This technique, in which the nondifferentiable 
double obstacle nonlinearity is approximated by differentiable logarithmic nonlinearities,
was first developed by P.~Colli, M.H.~Farshbaf-Shaker and J.~Sprekels in the paper
{\sl ``A deep quench approach to the optimal control of
an Allen--Cahn equation with dynamic boundary conditions and double obstacles''}
(Appl.\ Math.\ Optim.\ {\bf 71} (2015), pp.~1-24)
and has proved to be a powerful tool in a number of optimal control problems with
double obstacle potentials in the framework of systems of Cahn--Hilliard type. 
We first give a general convergence analysis of the deep quench approximation that includes an error estimate 
and then demonstrate that its use  leads in the double obstacle case to appropriate 
first-order necessary optimality
conditions in terms of a variational inequality and the associated adjoint state system.
    
\End{abstract}
\salta
\pagestyle{myheadings}
\newcommand\testopari{\sc Colli \ --- \ Gilardi \ --- \ Sprekels}
\newcommand\testodispari{\sc Deep quench approximation of  fractional Cahn--Hilliard systems}
\markboth{\testopari}{\testodispari}
\finqui
%

\section{Introduction}
\label{Intro}
\setcounter{equation}{0}

Let $\Omega\subset \erre^3$ denote an open, bounded, and connected set with smooth
boundary $\,\Gamma\,$  and outward normal derivative $\,\dn$, let $T>0$  be a 
final time,  and let $H:=L^2(\Omega)$ denote the Hilbert space of square-integrable 
real-valued functions defined on $\Omega$, endowed with the standard inner product
$(\cdot,\cdot)$ and norm $\,\|\cdot\|$, respectively. We denote 
\,$Q_t:=\Omega\times (0,t)\,$ for $0<t<T$ and \,$Q:=\Omega\times (0,T)$. 
We investigate in this paper the approximation and optimal control of an 
abstract system of evolutionary variational (in)equalities. More precisely, the variational
state system has the following form: we look 
for functions 
$(\mu,y)$ such that  
\begin{align}
\label{regy}
&y\in H^1(0,T;V_A^{-r})\cap L^\infty(0,T;V_B^\sigma) \quad\mbox{and} \quad \tau\dt y\in L^2(0,T;H)\\[0.5mm]
\label{regmu}
&\mu\in L^2(0,T;V_A^r),\\[0.5mm]
\label{L1}
&f_1(y)\in L^1(Q),
\end{align}
and \pier{satisfying}
\begin{align}
\label{weak1}
&\langle \dt y(t),v\rangle_{A,r}\,+\,(A^r\mu(t),A^r v)\,=\,0 \quad
\mbox{for every $\,v\in V_A^r\,$ and a.e. $t\in (0,T)$,} \\[1mm]
\noalign{\allowbreak}
\label{weak2}
&(\tau\dt y(t),y(t)-v)\,+\,(B^\sigma y(t),B^\sigma(y(t)-v))\,+\iO f_1(y(t))\nonumber\\
&\quad {}+\,(f_2'(y(t))-u(t), y(t)-v)\,\le\,(\mu(t),y(t)-v)\,+\iO f_1(v)\nonumber\\[1mm]
\noalign{\allowbreak}
&\mbox{for every \,$v\in V_B^\sigma\,$ and a.e. \,$t\in (0,T)$},\\[1mm]
\label{weak3}
&y(0)=y_0 \quad\mbox{in \,$\Omega$}.
\end{align}
Here, it is understood that 
$\,\iO f_1(v)=+\infty\quad\mbox{whenever}\quad f_1(v)\not\in\Luno.$ 
The precise meaning of the involved quantities and spaces will be given below. Notice that
\eqref{weak1}--\eqref{weak3} is a generalized version of the evolutionary system
\begin{align}
  & \dt y + A^{2r} \mu = 0\quad\mbox{in \,$Q$},
  \label{Iprima}
  \\
  & \tau \dt y + B^{2\sigma}  y + \pier{\partial f_1}(y)+f_2'(y) \pier{{}\ni{}} \mu + u\quad\mbox{in \,$Q$},
  \label{Iseconda}
  	\\
  & y(0) = y_0 \quad\mbox{in \,$\Omega$}.
  \label{Icauchy}
\end{align}
Here, $\tau\ge 0$ is a constant, $f_2:\erre\to\erre$ is a smooth function, and $f_1:\erre\to [0,+\infty]$ denotes a proper, convex, 
and lower semicontinuous function with $f_1(0)=0$, whose effective domain $D(f_1)$ is a closed interval in $\erre$ (possibly $\erre$ itself) 
and which is smooth in the interior of $D(f_1)$. \pier{In \eqref{Iseconda}, $\partial f_1$ denotes the subdifferential of $f_1$, which is a 
multivalued operator, in general, so that the inclusion replaces the equality.} 
The linear operators $A^{2r}$, and $B^{2\sigma}$, 
with $r>0$ and $\sigma>0$, 
denote fractional powers (in the spectral sense) of operators $A$ and $B$. We will 
give a proper definition of such operators in the next section. Throughout this 
paper, we generally assume: 
 
\vspace{1mm}\noindent
{\bf (A1)} \,\,\,$A:D(A)\subset H\to H$\, and \,$B:D(B)\subset H\to H\,$ are 
unbounded, monotone, and selfadjoint linear operators with compact resolvents.

\vspace{1mm}\noindent
This assumption implies that there are sequences 
$\{\lambda_j\}$ and $\{\lambda'_j\}$ of eigenvalues
and orthonormal sequences $\{e_j\}$ and $\{e'_j\}$ of corresponding eigenvectors,
that~is,
\Beq
  A e_j = \lambda_j e_j, \quad
  B e'_j = \lambda'_j e'_j,
  \aand
  (e_i,e_j) = (e'_i,e'_j) = \delta_{ij},
  \quad \hbox{for $i,j=1,2,\dots$,}
  \label{eigen}
\Eeq
\pier{with $\delta_{ij}$ denoting the Kronecker index,} such that
\begin{align}
  & 0 \leq \lambda_1 \leq \lambda_2 \leq \dots,
  \aand
  0 \leq \lambda'_1 \leq \lambda'_2 \leq \dots,
  \quad \hbox{with} \quad
  \lim_{j\to\infty} \lambda_j
  = \lim_{j\to\infty} \lambda'_j
  = + \infty,
  \label{eigenvalues}
  \\[1mm]
  & \hbox{$\{e_j\}$ and $\{e'_j\}$ are complete systems in $H$}.
  \label{complete}
\end{align}
  
The state system \eqref{Iprima}--\eqref{Icauchy} (and thus also \eqref{weak1}--\eqref{weak3})
can be seen as a generalization
of the famous Cahn--Hilliard system which models a  
phase separation process taking place in the container $\Omega$. In this case,
one typically has  $A^{2r}=B^{2\sigma}=-\Delta$ with zero Neumann or Dirichlet boundary conditions, and the unknown 
functions $\,y\,$ and $\,\mu\,$ stand for the \emph{order parameter} (usually a scaled
density of one of the involved phases) and the \emph{chemical potential}
associated with the phase transition, respectively. Moreover, $f:=f_1+f_2$ is a double-well potential. Typical cases
are the {\em classical regular potential}, the {\em logarithmic potential\/}, and the {\em double obstacle potential}, which (in this order)
are given by 
\begin{align}
\label{regpot}
& f_{{\rm reg}}(v) := \frac 14 \, (v^2-1)^2 \,,
  \quad v \in \erre, 
  \\[2mm]
\label{logpot}
& f_{{\rm log}}(v) := \left\{
  \begin{array}{ll}
  (1+v)\ln (1+v)+(1-v)\ln (1-v) - c_1 v^2 
  &\quad\mbox{for }\, v \in (-1,1)\\
  2\ln(2)-c_1&\quad\mbox{for }\,v\in\{-1,1\}\\
  +\infty&\quad\mbox{for }\,v\not\in [-1,1]
  \end{array}
  \right.
\\[2mm]
\label{obspot}
&f_{\rm obs}(v):= \left\{
\begin{array}{ll}
-c_1\,v^2 &\quad\mbox{if }\,|v|\le 1\\
+\infty &\quad\mbox{otherwise}
\end{array}
\right.
\end{align}
Here the constant $c_1>0$ is such that the above potentials are nonconvex.

Recently, in \cite[Thm.~2.6~and~2.8]{CGS18}, it was shown that the system  \eqref{weak1}--\eqref{weak3}  
 admits a solution $(\mu,y)$ satisfying \eqref{regy}--\eqref{L1}, where the admissible nonlinearities include all of the three cases 
\eqref{regpot}--\eqref{obspot}. In the analysis, it turned out that the first eigenvalue $\lambda_1$ of $\,A\,$ plays an important role.
 Indeed, the main assumption for the operators $A,B$  besides {\bf (A1)} was the following:

\vspace{1mm}\noindent
{\bf (A2)} \,\,\,Either\\
\hspace*{20mm}(i) \,\,\,$\lambda_1>0$\\
\hspace*{12.5mm} or\\
\hspace*{20mm}(ii) \,\,$0=\lambda_1<\lambda_2$, and $\,e_1\,$ is a constant and belongs to the domain of $\,B^\sigma$. 

\vspace{1mm}\noindent
The existence proof in \cite{CGS18} was based on Moreau--Yosida approximation, which is generally applicable to all of 
the three cases \eqref{regpot}--\eqref{obspot}. It turned out that the second component \,$y$\, of the solutions $(\mu,y)$ 
is always uniquely determined, while this is not necessarily so for the chemical potential $\mu$ (for cases in which also 
$\,\mu\,$ is unique, see \cite[Rem.~4.1]{CGS18} and \cite[Rem.~3.4]{CGS19}).   

In this paper, we focus on the case when $\,f=f_{{\rm obs}}$, that is, when $\,f_1=\indi\,$
is the indicator function of the interval $[-1,1]$, given by $\indi(v)=0$ if $v\in [-1,1]$ and 
$\indi(v)=+\infty$ otherwise. In this case, any solution $(\mu,y)$ of \eqref{weak1}--\eqref{weak3} must satisfy  \,$\int_Q \indi(y) <+\infty$, 
which entails that $y\in [-1,1]$ almost everywhere in $Q$ and thus  \,$\iO f_1(y(t))=0\,$ for almost every $t\in (0,T)$ in \eqref{weak2}.   
 While the question of well-posedness was settled in \cite[Thm.~2.6~and~2.8]{CGS18} for $\,f_1=\indi$, the matter of optimal control is still open. 
Indeed, the optimal control theory recently developed in \cite{CGS19} applies to certain classes of differentiable potentials only. In this paper, 
we aim at extending this theory to the case $\,f=f_{{\rm obs}}$. More precisely, we investigate  the following 
optimal control problem:

\vspace{1mm}
\noindent \CP0 \quad  Minimize the tracking-type cost functional
\begin{align}
\label{cost}
{\cal J}(y,u):=& \frac{\beta_1}2\,\|y(T)-y_\Omega\|^2\,+\,\frac {\beta_2}2\int_0^T\!\|y(t)-y_Q(t)\|^2\,dt
\,+\,\frac{\beta_3}2\int_0^T\!
\|u(t)\|^2\,dt
\end{align}
over the admissible set
\begin{equation}
\label{uad}
\uad:=\left\{u\in H^1(0,T;\Ldue):\,|u|\le\rho_1 \mbox{ \,a.\,e. in }\,Q, \,\,\,\,\|u\|_
{H^1(0,T;\Ldue)}\,\le\,\rho_2\right\},
\end{equation}
subject to \eqref{weak1}--\eqref{weak3} with $f_1=\indi$.
Here, $\rho_1>0$ and $\rho_2>0$ are such that $\uad\not= \emptyset$, $\beta_i$, $i=1,2,3$, are 
nonnegative but not all zero, and the given target functions satisfy $y_\Omega\in \pier{L^2(\Omega)} $ and $y_Q\in L^2(Q)$. Note
that \CP0 is well defined, since the component \,$y$\, of the solutions to the state system is uniquely determined.

The main difficulty inherent in \CP0 is the nondifferentiability of the nonlinearity $\indi$, which
entails that standard constraint qualifications from optimal control theory are violated, so that 
suitable Lagrange multipliers cannot easily be constructed. In such situations, the 
so-called ``deep quench'' approximation
has  proved to be a useful tool in a number of cases in the framework
of Cahn--Hilliard systems (see, e.g.,\cite{CFGS1, CFGS3, CGSEECT, CGSConvex, CS}). In all of these works, the starting point was that the optimal control problem (we will later
denote this problem by \CPal) had been successfully treated (by proving the Fr\'echet differentiability of the control-to-state operator and establishing
first-order necessary optimality conditions in terms of a variational inequality and the adjoint
state system) 
for the case when in the state system \eqref{weak1}--\eqref{weak3} the nonlinearity 
$f_1=\indi$ is \juerg{for $\alpha>0$} replaced by $f_1=h^\alpha:=\varphi(\alpha)h$, with the functions
\begin{eqnarray}\label{defphi}
\hspace*{-10.5mm}\mbox{\bf (A3)} \,&&\phi\in C^1[0,+\infty) \,\mbox{ is strictly increasing and satisfies \,} 
\lim_{\alpha\searrow0}\phi(\alpha)=0\pier{;} \hspace*{16mm}\\
\label{defh} 
&&h(v)=\left\{
\begin{array}{ll}
(1+v)\ln(1+v)+(1-\juerg{v})\ln(1-v), \,\,\,&\pier{v\in(-1,1)}\\
2\ln(2),\,\,\,&v\in\{-1,1\}\\
+\infty, \,\,\,&v\not\in [-1,1]
\end{array}
\right. \pier{.}
\end{eqnarray}
\juerg{We obviously have that}
\begin{align}
\label{monoal}
&0\,\le\,h^{\alpha_1}(v)\,\le\,h^{\alpha_2}(v)\quad\forall\,v\in\erre, \quad\mbox{if }\,0<\alpha_1<\alpha_2, \\[1mm]
\label{limh}
&\lim_{\alpha\searrow0} h^\alpha(v)\,=\,\indi(v)\quad\forall\,v\in \erre.
\end{align}
In addition, $\,h'(v)=\ln\left(\frac{1+v}{1-v}\right)$ \,and\, $h''(v)=\frac 2{1-v^2}>0$\, for
$v\in (-1,1)$, and thus, in particular,
\begin{align}
\label{limh1}
&\lim_{\alpha\searrow0}\,\phi(\alpha)h'(v)=0 \quad\mbox{for }\,-1<v<1,\\
\label{limh2}
&\lim_{\alpha\searrow0} \Bigl(\phi(\alpha)\,\lim_{v\searrow -1}h'(v)\Bigr)=-\infty, \quad
\lim_{\alpha\searrow0} \Bigl(\phi(\alpha)\,\lim_{v\nearrow +1}h'(v)\Bigr)=+\infty.
\end{align} 
We may therefore regard the graphs of the single-valued functions 
\begin{equation}
(h^\alpha)'(v)\,=\, \phi(\alpha)h'(v), \quad \mbox{for}\quad v\in (-1,1)\quad\mbox{and}\quad \alpha>0,
\end{equation} 
as approximations to the graph of the multi-valued subdifferential $\dI$.
Now the 
well-posedness results of  \cite{CGS18,CGS19} apply, 
yielding a solution pair $(\mual,\yal)$ for every $\alpha>0$, where the component $\yal$ is uniquely
determined. It is a natural 
question whether we have $y^\alpha\to y$ as $\alpha\searrow0$ in a suitable topology.
Below (cf.~Theorem 3.5), we will show that this is actually true; in Corollary 3.6, we will show that in a 
very special case with some global constant $K_2>0$ a quantitative  error estimate of the form
\begin{equation}
\|\yal-y \|_{C^0([0,T];\Ldue)\cap L^2(0,T;\Huno)}\,\le\,K_2\,|\alpha|^{1/2}
\end{equation}   
is valid. Also, owing to the construction, the \pier{approximate functions} $\,y^\alpha\,$  automatically attain values in
the domain of $\indi$; that is, we have $\|y^\alpha \|_{L^\infty(Q)}\,\le\,1$ \,for
all $\alpha>0$.

As far as the optimal control problem is concerned, the general strategy is then to derive 
uniform (with respect to $\alpha\in (0,1]$) a priori estimates for the state and adjoint state variables of 
\juerg{an ``adapted'' version of} \CPal\ that are sufficiently strong as to permit a passage to the limit as $\alpha\searrow0$ in order to derive meaningful first-order necessary optimality conditions also for \CP0. 
We can follow this strategy in this paper, since  in \cite{CGS19}  a 
corresponding theory for \CPal\ with the logarithmic potential \eqref{logpot} has been developed. 

However, while the approximation results for the solutions of the state system hold true under essentially the same 
assumptions as those imposed in \cite{CGS18} for the well-posedness results, it seems impossible  
to prove that the control-to-state operator  is Fr\'echet differentiable
between suitable Banach spaces without having at disposal suitable uniform $L^\infty(Q)$ bounds for both the 
state component $y$ and the functions $\,f^{(i)}(y)$, 
for $i=1,2,3$. In the case of the logarithmic potentials \,$h^\alpha$, which we intend to use for the
deep quench approximation, this means that we need to separate $\,\yal\,$ away
from the critical arguments $\pm 1$. Unfortunately, this postulate has the unpleasant consequence that 
Fr\'echet differentiability (and thus satisfactory first-order necessary optimality conditions) can only 
be established under rather restrictive conditions on the operators $A$ and $B$. A particular 
case in which our analysis will work is given if $A=B=-\Delta$ with zero Neumann boundary condition, 
$\sigma=1/2$, and $r\ge 3/8$.  

Let us add a few remarks on the existing literature. \pier{One can find} numerous contributions on viscous/nonviscous,
local/nonlocal, convective/nonconvective Cahn--Hilliard systems for the classical (non-fractional) case 
$A=B=-\Delta$, $2r=2\sigma=1$, where various types of boundary
conditions (e.g., Dirichlet, Neumann, dynamic) and different assumptions on the nonlinearity  were
considered. We refer the
interested reader to the recent paper \cite{CGSAnnali} for a selection of associated references. Some papers
also address the coupled Cahn--Hilliard/Navier--Stokes system (see, e.g, \cite{FGG}, \cite{FGGS}, and the 
references given therein).

The literature on optimal control problems for non-fractional Cahn--Hilliard systems is still rather scarce.    
The case of Dirichlet and/or Neumann boundary conditions for various types of such systems were the
subject of, e.g., the works\cite{CGRS, CGSAIMS, CGSAMO, CGSEECT, Duan, WN, Z, ZW}, while the case of  
dynamic boundary conditions was studied in\cite{CFGS1, CFGS2, CFGS3, CGSAdvan, CGSAMO, CGSSIAM18, CGSConvex,
CS, Fukao}. The optimal control of convective Cahn--Hilliard systems was addressed in\cite{CGSSIAM18,
CGSConvex,RS, ZL1, ZL2},
while the papers\cite{BDM1, BDM2, FGS, FRS, HHCK, HKW, HW1, HW2, HW3, Medjo} were concerned with 
coupled Cahn--Hilliard/Navier--Stokes systems. 

There are only a few contributions to the theory of Cahn--Hilliard systems involving fractional operators.
In the connection of well-posedness and regularity results, we refer to \cite{AM,AkSeSchi} for the case of 
the fractional negative Laplacian with zero Dirichlet boundary conditions; general operators other than the 
negative Laplacian
have apparently only been 
studied in\gianni{\cite{CGS18,GalDCDS,GalEJAM,GalAIHP}}. 
As of now, aspects of optimal control have been scarcely
dealt with even for simpler linear evolutionary systems involving fractional operators; for such systems,
some identification problems were addressed in the recent contributions \cite{GV, SV}, while
for optimal control problems for such cases we refer to~\cite{HAS} (for the stationary (elliptic) case, see
also \juerg{\cite{AKW,HO1, HO2,APR,AW1,AW2})}.
However, to the authors' best
knowledge, the present paper appears to be the first contribution that addresses optimal control problems
for Cahn--Hilliard systems with general fractional order operators and potentials of double obstacle
type.

The paper is organized as follows: the subsequent Section~2 brings some auxiliary functional analytic 
material on fractional order operators, while in Section~3  we establish some general 
convergence results for the deep quench 
approximation of the state system \juerg{\eqref{weak1}--\eqref{weak3}}. In particular, an error
estimate is proved. 
In Section~4, we investigate the relations between the solutions to the optimal control problems
(${\mathcal{P}}_0$) and the solutions to the corresponding optimal control problems for the deep
quench approximations. In the final Section~5, we then employ the results from \cite{CGS19}
to establish the first-order necessary optimality conditions 
for (${\mathcal P}_0$).

Throughout the paper, we denote for a general Banach space $\,X\,$ other than $H=\Ldue$ by 
$\,\|\cdot\|_X\,$ and  $\,X^*\,$ its norm 
and dual space, respectively; the dual pairing between elements of $X^*$ and $X$ is denoted by $\langle\cdot,\cdot
\rangle_X$. 


\section{Fractional powers and auxiliary results}
\label{FPAM}
\setcounter{equation}{0}

In this section, we collect some auxiliary material concerning functional analytic notions.
To this end, we generally assume that the conditions {\bf (A1)} and {\bf (A2)} are satisfied.
At this point, some remarks on the assumption {\bf (A2)} are in order.

\Brem
The condition $\lambda_1>0$  is satisfied
for many standard elliptic operators of second or higher order with zero Dirichlet boundary conditions
(however, also zero mixed boundary conditions could be considered, 
with proper definitions of the domains of the operators); typical cases are 
the (negative) Laplacian \,$A=-\Delta$ with the domain $D(-\Delta)=\Hdue\cap\Hunoz$ or
the bi-harmonic operator $A=\Delta^2$ with the domain\, 
$D(\Delta^2)=\Hx4\cap H^2_0(\Omega)$. On the other hand, we have $0=\lambda_1<\lambda_2$ and
$e_1\equiv {\rm const.}$ \,for important problems with zero Neumann boundary conditions;
typical examples are \,$A=-\Delta$\, with the domain $D(-\Delta)=\{v\in\Hdue:\ \dn v=0 \mbox{ \,on\, }\Gamma\}$\, and \,$\,A=\Delta^2$ with the domain  
$D(\Delta^2)=\{v\in\Hx4:\ \dn v=\dn\Delta v=0\mbox{ \,on\, }\Gamma\}$. We also point out that
$A$ and $B$ can be completely unrelated if $\lambda_1>0$, while in the other case
the constant functions have to belong to $D(B^\sigma)$. The latter holds true if $B=-\Delta$
with the domain  $D(-\Delta)=\{v\in\Hdue:\ \dn v=0 \mbox{ \,on\, }\Gamma\}$, while in the Dirichlet case
$D(-\Delta)=\Hdue\cap H_0^1(\Omega)$ no nontrivial constant functions are contained in
$D(B)$; however, \gianni{if $0<\sigma<1/4$, \juerg{then} $D(B^\sigma)$~coincides with the usual Sobolev--Slobodeckij space $\Hx{2\sigma}$
and thus contains all constant functions}.
\Erem

Using the facts summarized in \eqref{eigen}--\eqref{complete},
we can define the powers of $A$ and $B$ 
for an arbitrary positive real exponent.
For the first operator, we have
\Bsist
  && \VA r := D(A^r)
  = \Bigl\{ v\in H:\ \somma j1\infty |\lambda_j^r (v,e_j)|^2 < +\infty \Bigr\}
  \aand
  \label{defdomAr}
  \\
  && A^r v = \somma j1\infty \lambda_j^r (v,e_j) e_j
  \quad \hbox{for $v\in\VA r$},
  \label{defAr}
\Esist
the series being convergent in the strong topology of~$H$,
due to the properties \eqref{defdomAr} of the coefficients.
In principle, we can endow $\VA r$ with the (graph) norm and inner product
\Beq
  \norma v_{gr,A,r}^2 := (v,v)_{gr,A,r}
  \aand
  (v,w)_{gr,A,r} := (v,w) + (A^r v , A^r w)
  \quad \hbox{for $v,w\in\VA r$}.
  \label{defnormagrAr}
\Eeq
This makes $\VA r$ a Hilbert space.
However, we can choose any equivalent Hilbert norm.
Indeed, in view of assumption {\bf (A2)}, it is more convenient to work with the Hilbert norm 
\Beq
  \norma v_{A,r}^2 := \left\{ 
  \begin{aligned}
  & \norma{A^r v}^2
  = \somma j1\infty |\lambda_j^r (v,e_j)|^2
  \qquad \hbox{if $\lambda_1>0$,}
  \\
  & |(v,e_1)|^2 + \norma{A^r v}^2
  = |(v,e_1)|^2 + \somma j2\infty |\lambda_j^r (v,e_j)|^2
  \qquad \hbox{if $\lambda_1=0$}.
  \end{aligned}
  \right.
  \label{defnormaAr}
\Eeq
In \cite[Prop.~3.1]{CGS18} it has been shown that this norm is equivalent 
to the graph norm defined in~\eqref{defnormagrAr},
and we always will work with the norm \eqref{defnormaAr} instead of with \eqref{defnormagrAr}.
We also use the corresponding inner product in $\VA r$ 
given~by
\Bsist
  && (v,w)_{A,r}
  = (A^r v,A^r w)
  \quad \hbox{or} \quad
  (v,w)_{A,r}
  = (v,e_1)(w,e_1) + (A^r v,A^r w),
  \qquad
  \non
  \\
  \label{inpro}
  &&  \hbox{depending on whether $\lambda_1>0$ or $\lambda_1=0$,\quad 
	for $v,w\in\VA r$.}
\Esist

\Brem
Observe that in the case $\lambda_1=0$ 
the constant value of $e_1$
equals one of the numbers $\pm|\Omega|^{-1/2}$,
where $|\Omega|$ is the volume of~$\Omega$.
It follows for every $v\in H$ that the first term $(v,e_1)e_1$ of the Fourier series of $v$
is the constant function whose value~is the mean value of~$v$, which is defined by
\Beq
  \mean (v) := \frac 1 {|\Omega|} \iO v\,.
  \label{defmean}
\Eeq
\Erem

In the same way as for $A$, starting from \accorpa{eigen}{complete} for~$B$,
we can define the power $B^\sigma$ of $B$ for every $\sigma>0$, where for 
$V_B^\sigma$ we choose the graph norm.
We therefore set
\Bsist
  && \VB\sigma := D(B^\sigma),
  \quad \hbox{with the norm $\norma\cpto_{B,\sigma}$ associated to the inner product}
  \label{defBs}
  \non
  \\
  && (v,w)_{B,\sigma} := (v,w) + (B^\sigma v,B^\sigma w)
  \quad \hbox{for $v,w\in \VB\sigma$}.
  \label{defprodBs}
\Esist

To resume our preparations, we observe that if $r_i$ and $\sigma_i$ are arbitrary positive exponents,
then it is easily seen that we have the ``Green type'' formulas
\Bsist
  && (A^{r_1+r_2} v,w)
  = (A^{r_1} v, A^{r_2} w)
  \quad \hbox{for every $v\in\VA{r_1+r_2}$ and $w\in\VA{r_2}$},
  \label{propA}
  \\[1mm]
  && (B^{\sigma_1+\sigma_2} v,w)
  = (B^{\sigma_1} v, B^{\sigma_2} w)
  \quad \hbox{for every  $v\in\VB{\sigma_1+\sigma_2}$ and $w\in\VB{\sigma_2}$}.
  \label{propB}
\Esist

The next step is the introduction of some spaces with negative exponents.
We set
\Beq
  \VA{-r} := (\VA r)^* 
  \quad \hbox{for $r>0$},
  \label{defVAneg}
\Eeq
and endow $\VA{-r}$ with the dual norm $\,\|\cdot\|_{A,-r}\,$ of $\,\|\cdot\|_{A,r}$.
We use the symbol $\langle\cpto,\cpto\rangle_{A,r}$ for the duality pairing
between $\VA{-r}$ and~$\VA r$ and identify $H$ with a subspace of $\VA{-r}$
in the usual way, i.e., such that
$\,\langle v,w \rangle_{A,r} = (v,w)\,$ for every $v\in H$ and $w\in\VA r$. Likewise,
we~set
\Beq 
V_B^{-\sigma}:=(V_B^\sigma)^* \quad\mbox{for }\,\sigma>0.
\Eeq
As $\,V_B^\sigma\,$ is dense in $\,H$, we have the analogous embedding
\Beq
H\subset V_B^{-\sigma}.
\Eeq

Observe that the following embedding results are valid:
\begin{align}
  & \hbox{The embeddings $\VA{r_2} \subset \VA{r_1} \subset H$ are dense and compact for $0<r_1<r_2$}.
  \label{compembA}
  \\
  & \hbox{The embeddings $H \subset \VA{-r_1} \subset \VA{-r_2}$ are dense and compact for $0<r_1<r_2$}.
  \qquad
  \label{compembAneg}
  \\
  & \hbox{The embeddings $\VB{\sigma_2} \subset \VB{\sigma_1} \subset H$ are dense and compact for $0<\sigma_1<\sigma_2$}.
  \label{compembB}
\end{align}

At this point, we introduce the Riesz isomorphism $\calR_r:\VA r\to\VA{-r}$
associated with the inner product~\eqref{inpro}, which is given by
\Beq
  \< \calR_r v , w >_{A,r}
  = (v,w)_{A,r}
  \quad \hbox{for every $v,w\in\VA r$}.
  \label{riesz}
\Eeq
Moreover, we~set
\begin{align}
  & \Vz r := \VA r
  \aand
  \Vz{-r} := \VA{-r}
  \quad \hbox{if $\lambda_1>0$},
  \non
  \\[1mm]
  & \Vz r := \{v\in \VA r :\ \mean (v)=0\}
  \aand
  \Vz{-r} := \{v \in \VA{-r} :\ \<v,1>_{A,r}=0 \}
  \quad \hbox{if $\lambda_1=0$} \,.
  \label{defVrpos}
\end{align}
According to \cite[Prop.~3.2]{CGS18},  $\calR_r$ maps $\Vz r$ onto $\Vz{-r}$
and extends to $\Vz r$ the restriction of $A^{2r}$ to $\Vz{2r}$.
In view of this result, it is reasonable to use a proper notation
for the restrictions of $\calR_r$ and $\calR_r^{-1}$ to the subspaces
$\Vz r$ and $\Vz{-r}$, respectively.
We~set 
\Beq
  \Az{2r} := (\calR_r)_{|\Vz r}
  \aand
  \Az{-2r} := (\calR_r^{-1})_{|\Vz{-r}}\,,
  \label{defAz}
\Eeq
where the index $0$ has no meaning if $\lambda_1>0$ (since then $\Vz{\pm r}=\VA{\pm r}$),
while it reflects the zero mean value condition in the case $\lambda_1=0$.
We thus have
\Bsist
  && \Az{2r} \in \calL(\Vz r,\Vz{-r}) , \quad
  \Az{-2r} \in \calL(\Vz{-r},\Vz r)
  \aand
  \Az{-2r} = (\Az{2r})^{-1}\,,
  \label{contlinAz}
  \\[1mm]
  && \< \Az{2r} v,w >_{A,r} = (v,w)_{A,r} = (A^r v,A^r w)
  \quad \hbox{for every $v\in\Vz r$ and $w\in\VA r$}\,,
  \qquad
  \label{identityAz}
  \\[1mm]
  && \< f , \Az{-2r} f >_{A,r}
  = \norma{\Az{-2r} f}_{A,r}^2
  = \norma f_{A,-r}^2
  \quad \hbox{for every $f\in\Vz{-r}$}.
  \label{normaAz}
\Esist

\section{Deep quench approximation of the state system}

\setcounter{equation}{0}
In this section, we state our general assumptions and discuss the deep quench approximation
of the state system \juerg{\eqref{weak1}--\eqref{weak3}}.
Besides \gianni{{\bf (A1)}--{\bf (A3)}}, we generally assume for \gianni{the structure and} the data of the state system:

\vspace{2mm}\noindent
{\bf (A4)}  \,\,$r>0$, $\sigma>0$, and $\tau\ge 0$ are fixed real numbers.

\vspace{2mm}\noindent
{\bf (A5)} \,\,$f_2\in C^3(\erre)$, and  $\,f_2'\,$ is Lipschitz continuous on $\erre$ with Lipschitz constant
$L>0$.

\vspace{2mm}\noindent
{\bf (A6)} \,\,\pier{$y_0\in V^{\sigma}_B$}, \,  and \, $-1<\infess_{x\in\Omega}y_0(x),
\quad \supess_{x\in\Omega}y_0(x)<+1$.

\vspace{2mm}\noindent
{\bf (A7)} \,\,$u\in {\cal X}\,:= \,H^1(0,T;\Ldue)\cap L^\infty(Q)$.

\vspace{2mm}\noindent 
We draw a few consequences from {\bf (A6)}. Namely, the mean value of $y_0$ belongs to the interior
of both $D(\dI)$ and $D((h^\alpha)')$, for all $\alpha>0$. Moreover, we have $\indi(y_0)\in\Luno$ and
$h(y_0)\in\Luno$, and
$h'(y_0)$ belongs to $\Ldue$. Thus, the conditions \cite[(2.27), (2.28)]{CGS18} 
on $y_0$ for the application of \cite[Thm.~2.6]{CGS18}
are satisfied, where we note that 
\begin{equation}
\label{estini}
\|h^\alpha(y_0)\|_{\Luno}\,+\,\|(h^\alpha)'(y_0)\|_{\Ldue}\,\le\,\hat c \quad\forall\,\alpha\in (0,1],
\end{equation}
with some constant $\hat c>0$ which is independent of $\alpha\in (0,1]$.

\vspace{2mm}
We now consider the state system \eqref{weak1}--\eqref{weak3} for the cases $f_1=\indi$ and $f_1=h^\alpha$ ($\alpha
\in (0,1]$), respectively. By virtue of \cite[Thm.~2.6]{CGS18}, there exist solution pairs \,$(\mu,y)$\, and 
$\,(\mual,\yal)$, respectively, which enjoy the properties \eqref{regy}--\eqref{L1}, and the
(uniquely determined) second components satisfy 
\begin{equation}
\label{interv}
-1\le y\le 1\, \mbox{ a.e. in $\,Q$,}\quad -1\le \yal\le 1\,\mbox{ a.e. in \,$Q$.}
\end{equation} 
We are now going to investigate the behavior of the family
$\{(\mual,\yal)\}_{\alpha>0}$ of deep quench approximations for $\alpha\searrow0$. We begin our analysis with the derivation of
 general a priori estimates.

\Bthm
Suppose that the general assumptions {\bf (A1)}--{\bf \gianni{(A7)}} are
fulfilled, and assume that $(\mual,\yal)$ are solution pairs to the problem \eqref{weak1}--\eqref{weak3} with
$f_1=h^\alpha$ for $\alpha\in (0,1]$ as established in {\rm \cite[Thm.~2.6]{CGS18}}.
Then there exists a constant $K_1>0$, which only depends on the data of the system
\eqref{weak1}--\eqref{weak3}, such that 
\begin{align}
\label{albound1}
&\|\mual\|_{L^2(0,T;V_A^{r})}\,+\,\|\yal\|_{H^1(0,T;V_A^{-r})\cap L^\infty(0,T;V_B^\sigma)\cap 
L^\infty(Q)}
\,+\,\|\phi(\alpha)h(\yal)\|_{L^\infty(0,T;\Luno)}\nonumber\\
&+\,\|\tau^{1/2}\dt\yal\|_{L^2(0,T;H)}\,\le\,K_1\quad\forall\,\alpha\in (0,1].
\end{align}
If, in addition, 
\begin{equation}
\label{adco}
\tau>0  \quad\mbox{and}\quad y_0\in V_B^{2\sigma}, 
\end{equation}
then we have the additional bounds
\begin{align}
\label{albound2}
&\|\yal\|_{W^{1,\infty}(0,T;H)\cap H^1(0,T;V_B^\sigma)}\,+\,\|\mual\|_{L^\infty(0,T;V_A^{2r})}
\nonumber\\
&+\int_0^T\!\!\iO\phi(\alpha)h''(\yal)\left|\dt\yal\right|^2\,\le\,K_1
\quad\forall\,\alpha\in (0,1]. 
\end{align}
\Ethm

\Bdim
To establish the validity of \eqref{albound1}, we have to follow the lines of the proof of \cite[Thm.~2.6]{CGS18}.
The method of proof of \cite[Thm.~2.6]{CGS18}, specified to our situation where the convex
part of the nonlinearity is given by $h^\alpha$, was the following:

\noindent {\sc Step 1:} Replace in \eqref{weak2} the function $\,f_1=h^\alpha$\, by its Moreau--Yosida
approximation~$\,h^\alpha_\lambda$, where $\lambda>0$.\\[0.5mm]
{\sc Step 2:} Approximate the resulting system of variational inequalities (which on the level of 
Moreau--Yosida approximations become variational equalities) via time discretization.\\[0.5mm]
{\sc Step 3:} Show unique solvability for the discrete system and derive a priori estimates for
the discrete approximations.\\[0.5mm]
{\sc Step 4:} Take the time step-size to zero in the time-discrete system to establish 
unique solvability of the system governing the Moreau--Yosida approximations.\\[0.5mm]
{\sc Step 5:} Take the limit as $\lambda\searrow0$ to obtain the solvability of the system \eqref{weak1}--\eqref{weak3}
for $f_1=h^\alpha$.
\\[0.5mm]
{\sc Step 6:} Show the uniqueness of the second solution component $\yal$.

Now, a closer inspection reveals that in our case all of the bounds established in the a priori estimates 
performed in {\sc Step 3} are uniform with respect to $\alpha\in (0,1]$, and due to the semicontinuity properties
of norms, they persist under the limit processes as the step-size of the time discretization and
the Moreau--Yosida parameter $\lambda$ approach zero. The validity of the estimate \eqref{albound1} is thus a consequence
of the estimate \cite[Eq.~(6.1)]{CGS18}. 

To \pier{offer to} the reader a little flavor of the argument, we give a formal derivation of 
a part of \eqref{albound1} (which
becomes rigorous on the level of the time-discrete approximation). To this end, let us assume that $\dt\yal\in
L^2(0,T;V_B^\sigma)$\, (which is satisfied under the assumption \eqref{adco}) and that the variational inequality 
\eqref{weak2} with $f_1=h^\alpha$ is equivalent to the variational equation
\begin{align}
\label{al2new}
&(\tau\dt \yal(t),v)\,+\,(B^\sigma \yal(t),B^\sigma v)\,+\,((h^\alpha)'(\yal(t))+f_2'(\yal(t)),v)\,=\,
(\mual(t)+u(t),v)\nonumber\\[1mm]
&\mbox{for every \,$v\in V_B^\sigma\,$ and a.e. \,$t\in (0,T)$}.
\end{align}
The latter is certainly satisfied on the level of the Moreau--Yosida approximations to the deep
quench approximations $(\mual,\yal)$.  We then insert $v=\mual(t)$ in \eqref{weak1} (written for
$(\mu,y)=(\mual,\yal)$) and $v=\dt\yal(t)$ in \eqref{al2new},
add the resulting equations, and integrate with respect to time over $[0,t]$, where $t\in (0,T]$ is arbitrary.
It then follows after an obvious cancellation of terms that
\begin{align*}
&\tau\int_0^t\!\!\iO|\dt\yal|^2\,+\,\frac 12\,\|B^\sigma\yal(t)\|^2\,+\int_0^t\!\|A^r\mual(s)\|^2\,ds
\,+\iO h^\alpha(\yal(t))\\
&=\,\frac 12\,\|B^\sigma y_0\|^2\,+\iO(h^\alpha(y_0)-f_2(\yal(t))
+f_2(y_0))\,+\int_0^t\!\!\iO u\,\yal.
\end{align*}
Now we recall \eqref{estini} and the fact that $y_0\in V_B^\sigma$ (cf.~{\bf \gianni{(A6)}}). 
We thus can infer from \eqref{interv}, {\bf \gianni{(A5)}}, and {\bf \gianni{(A7)}}, 
that all of the terms on the \rhs\ are bounded independently of $\alpha\in (0,1]$
by a constant that depends in a continuous and monotone way on $\|u\|_{L^1(Q)}$. But this means
that
\begin{align*}
&\|\pier{A^r\mual}\|_{L^2(0,T;H)}+\|\yal\|_{L^\infty(0,T;V_B^\sigma)\cap L^\infty(Q)} 
+\|h^\alpha(\yal)\|_{\juerg{L^\infty(0,T;}\Luno)}\\
&+\|\tau^{1/2}\dt\yal\|_{L^2(0,T;H)}\,\le\,
C,
\end{align*}
where $C>0$ is independent of $\alpha\in (0,1]$. This is already a part of the asserted bound
\eqref{albound1}. Now, if $\lambda_1>0$, then \eqref{defnormaAr} and the above estimate 
immediately entail that $\{\mual\}_{\alpha\in(0,1]}$ is bounded in $L^2(0,T;V_A^r)$, and
comparison in \eqref{Iprima} yields a uniform bound for $\{\dt\yal\}_{\alpha\in (0,1]}$
in $L^2(0,T;V_A^{-r})$, which then shows that \eqref{albound1} is valid. In the case
when $\lambda_1=0$, the boundednes of  $\{\mual\}_{\alpha\in(0,1]}$ in $L^2(0,T;V_A^r)$
is shown by proving that the mean values of $\{\mual(t)\}_{\alpha\in(0,1]}$ are
uniformly bounded in $L^2(0,T)$.
For this argument, we refer the reader to the proof of \cite[Thm.~2.6]{CGS18}.

Assume now that also the condition \eqref{adco} is fulfilled. 
In order to prove the bounds \eqref{albound2}, we follow the proof
of \cite[Thm.~2.8]{CGS18}, which again uses the time-discrete approximation scheme for the
system governing the Moreau--Yosida approximations mentioned above in describing {\sc Step 3} in the proof
of \cite[Thm.~2.6]{CGS18}. At this point, we recall the estimate \eqref{estini}. With this estimate in mind,
it turns out that all of the estimates performed in the proof
of \cite[Thm.~2.8]{CGS18} on the discrete approximations yield bounds that do not depend on $\alpha\in (0,1]$
and persist under the limit processes of taking the time step-size and the Moreau--Yosida parameter
$\lambda$ to zero. Since \eqref{albound2} exactly reflects the bounds established there,
the assertion is proved. 

For the reader's convenience, we again provide a formal sketch of the argument. To this end, we formally
differentiate \eqref{al2new} with respect to
$t$, obtaining the identity
\begin{eqnarray}
\label{diffi1}
&&(\tau\partial^2_{tt}\yal,v)+(B^\sigma\dt\yal,B^\sigma v)+ (\phi(\alpha)h''(\yal)\dt\yal+f_2''(\yal)\dt\yal,v)
=(\dt\mual+\dt u,v)\non\\
&&\mbox{for every $\,v\in V_B^\sigma\,$ and a.e. in $\,(0,T)$}. 
\end{eqnarray}
Then we formally test \eqref{weak1} by $v=\dt\mual$ and \eqref{diffi1} by $v=\dt\yal$, and add the
resulting identities. After an obvious cancellation \juerg{of terms}, we arrive at 
\begin{align}
\label{diffi2}
&\frac \tau 2\,\|\dt\yal(t)\|^2\,+\,\frac 12\,\|A^r\mual(t)\|^2\,+\int_0^t\!\!\iO |B^\sigma\dt\yal|^2
\,+\int_0^t\!\!\phi(\alpha)h''(\yal)|\dt\yal|^2 \nonumber\\
&=\,\frac \tau 2\,\|\dt\yal(0)\|^2\,+\,\frac 12\,\|A^r\mual(0)\|^2
\,-\int_0^t\!\!f_2''(\yal)|\dt\yal|^2\,+\int_0^t\!\!\iO \dt u\,\dt\yal\,,
\end{align}
where the last summand on the left-hand side is nonnegative and the last two terms on the right-hand
side can be estimated by an expression of the form
$$C_1\int_0^t\!\!\iO\left(|\dt u|^2\,+\,|\dt\yal|^2\right)\,,$$
where $C_1>0$ is independent of $\alpha\in (0,1]$. We thus are left to estimate the initial value terms.
To this end, we formally write \eqref{weak1} and \eqref{al2new} for $t=0$, obtaining the
identities
\begin{align}
\label{diffi3}
&\langle\dt\yal(0),v\rangle_{A,r}+(A^r\mual(0),A^r v)=0 \quad\forall\,v\in V_A^r,\\[0.5mm]
\label{diffi4}
&(\tau\dt\yal(0),v)+(B^{2\sigma} y_0+(h^\alpha)'(y_0)+f_2'(y_0),v)\,=\,(\mual(0)+u(0),v)
\quad\forall\,v\in V_B^\sigma.
\end{align}
Now observe that, by virtue of \eqref{adco}, \eqref{estini}, and {\bf \gianni{(A5)}}, 
the sum $\,B^{2\sigma} y_0+ (h^\alpha)'(y_0)+f_2'(y_0)\,$ is bounded in $\Ldue$,
uniformly with respect to $\alpha\in (0,1]$. 
Hence, if we (formally) test \eqref{diffi3} by $\mual(0)$ and \eqref{diffi4} by $\dt\yal(0)$, add
the resulting identities, and apply Young's inequality (note that we have $\tau>0$ by 
assumption \eqref{adco}), then we arrive at an estimate of the form   
$$
\|A^r\mual(0)\|^2\,+\,\frac\tau 2\,\|\dt\yal(0)\|^2\,\le\, C_2\,\tau^{-1}(1+\|u(0)\|^2),
$$ 
where $C_2>0$ is independent of $\alpha\in (0,1]$. We may then combine this estimate with
\eqref{diffi2} to conclude from Gronwall's lemma that
\begin{align}
\|\dt\yal\|_{L^\infty(0,T;H)\cap L^2(0,T;V_B^\sigma)}\,+\,\|A^r\mual\|_{L^\infty(0,T;H)}
\,+\int_0^T\!\!\iO\phi(\alpha)h''(\yal)\left|\dt\yal\right|^2\,\le\,C_3,
\end{align}
where $C_3>0$ is independent of $\alpha\in (0,1]$. With this, the first and third summands on 
the left of \eqref{albound2} are uniformly bounded, which then, by comparison in \eqref{Iprima},
also holds true for $\,\|A^{2r}\mual\|_{L^\infty(0,T;H)}$. Hence, \eqref{albound2}
is proved if $\lambda_1>0$. In the case $\lambda_1=0$, it is necessary to derive
a uniform $L^\infty(0,T)$ bound for the mean values of $\,\{\mual(t)\}_{\alpha\in (0,1]}$.
\pier{About} this, we again refer the reader to the proof of \cite[Thm.~2.8]{CGS18}.	
\Edim

\Brem
A closer inspection of the a priori estimates for the time-discretized 
systems mentioned above \juerg{reveals} that the constant $K_1$ depends in a monotone and continuous  way
on the norm $\|u\|_{\cal X}$. Hence, for any bounded subset
${\cal U}$ of ${\cal X}$ (in particular, for ${\cal U}=\uad$) it follows that there is a constant, which is
again denoted by $K_1$, such that the estimates \eqref{albound1} and \eqref{albound2}, respectively, hold true whenever $u$ is an arbitrary element of ${\cal U}$.
\Erem
  
\vspace{2mm} Next, we show the convergence of the deep quench approximations. Before formulating the result, we notice
that the following control-to-state operators are well defined on the space ${\cal X}$:
\begin{align}
\label{defcs0}
\S_0&:\,{\cal X}\ni u\mapsto \S_0(u):= y,\\
\label{defcsal}
\S_{\alpha}&:\,{\cal X}\ni u\mapsto \S_{\alpha}(u):=\yal,
\end{align}
where $(\mu,y)$ and $(\mual,\yal)$ denote solutions to the systems \eqref{weak1}--\eqref{weak3} for $f_1=\indi$ and
$f_1=h^\alpha$, $\alpha\in (0,1]$, respectively, as established in \cite[Thm.~2.6]{CGS18}. We have the following result.

\Bthm
Suppose that the assumptions {\bf (A1)}--{\bf \gianni{(A7)}} are fulfilled, and let sequences $\{\alpha_n\}\subset (0,1]$ and 
$\{\un\}\subset{\cal X}$ be given such that $\alpha_n\searrow0$ and $\un\to u$ weakly-star in ${\cal X}$  as $n\to\infty$ 
for some $u\in{\cal X}$. Moreover, let $(\mn,\yn)$  
be solutions to \eqref{weak1}--\eqref{weak3} for $f_1=h^{\alpha_n}$ and $u=u_n$, $n\in\enne$, as established in 
{\rm \cite[Thm.~2.6]{CGS18}}. Then there are a solution $(\mu,y)$ with $y=\S_0(u)$ to the problem
\eqref{weak1}--\eqref{weak3} with $f_1=\indi$ and a subsequence $\{\alpha_{n_k}\}_{k\in\enne}$ of
$\{\alpha_n\}$ such that, 
as $k\to\infty$, 
\begin{align}
\label{conmu}
\mu^{\alpha_{n_k}}\to\mu&\quad\mbox{weakly in }\,L^2(0,T;V_A^r),\\
\label{cony}
y^{\alpha_{n_k}}\to y&\quad\mbox{weakly-star in }\,H^1(0,T;V_A^{-r})\cap L^\infty(0,T;V_B^\sigma)\nonumber\\
&\quad\mbox{and strongly in }\,C^0([0,T];H), \\
\label{conyt}
\dt y^{\alpha_{n_k}} \to \dt y&\quad\mbox{weakly in }\,L^2(0,T;H) \quad\mbox{if }\,\tau>0.
\end{align}
Moreover, if \eqref{adco} is fulfilled, then the above solution $(\mu,y)$ also satisfies
\begin{align}
\label{conmup}
\mu^{\alpha_{n_k}}\to\mu&\quad\mbox{weakly-star in }\,L^\infty(0,T;V_A^{2r}),\\
\label{conyp} 
y^{\alpha_{n_k}}\to y&\quad\mbox{weakly-star in }\,W^{1,\infty}(0,T;H)\cap H^1(0,T;V_B^\sigma).
\end{align} 
\Ethm
\Bdim
The sequence $\{\un\}$ converges weakly-star in ${\cal X}$ and thus forms a bounded subset of ${\cal X}$.
According to Remark~3.2, the bounds \eqref{albound1} and \eqref{albound2} (the latter if \eqref{adco}
is satisfied) apply, where the constant $K_1$ is independent of $n$. Therefore, there are limits $(\bm,\by)$  and a subsequence of $\{(\mn,\yn)\}$,
which is for convenience again indexed by $n$, such that, as~$n\to\infty$, 
\begin{align}
\label{p1}
\mn\to\bm&\quad\mbox{weakly in }\,L^2(0,T;V_A^r),\\
\label{p2}
\yn\to\by&\quad\mbox{weakly-star in }\,H^1(0,T;V_A^{-r})\cap L^\infty(0,T;V_B^\sigma),\\
\label{p2bis}
\pier{\yn\to\by}
&\quad\mbox{strongly in }\,C^0([0,T];H) \, \mbox{ and pointwise a.e. in }\,Q,\\
\label{p3}
\dt \yn \to \dt \by&\quad\mbox{weakly in }\,L^2(0,T;H) \quad\mbox{if }\,\tau>0,
\end{align}
and, if \eqref{adco} is satisfied, 
\begin{align}
\label{p4}
\mn\to\bm&\quad\mbox{weakly-star in }\,L^\infty(0,T;V_A^{2r}),\\
\label{p5} 
\yn\to \by&\quad\mbox{weakly-star in }\,W^{1,\infty}(0,T;H)\cap H^1(0,T;V_B^\sigma).
\end{align} 
Notice that the strong convergence result in \eqref{p2bis} follows from \cite[Sect.~8,~Cor.~4]{Simon}, since,
by \eqref{compembB}, $V_B^\sigma$ is compactly emdedded in $H$; we thus may without loss of generality 
assume that $\yn\to\by$ pointwise a.e. in $Q$. Since, by virtue of \eqref{interv}, $-1\le\yn\le+1$ a.e. in $Q$, 
we infer that $-1\le\by\le 1$ a.e. in $Q$, and thus  $\indi(\by)\in L^1(Q)$ with 
$$\int_0^T\!\!\int_\Omega\indi(\by)=0.$$

It remains to show that $(\bm,\by)$ is a solution to \eqref{weak1}--\eqref{weak3} in the sense of 
\cite[Thm.~2.6]{CGS18} for $f_1=\indi$ and control $\,u$. To this end, we pass to the limit as $n\to\infty$ in the system \eqref{weak1}--\eqref{weak3},
written for $f_1=h^{\alpha_n}$ and $u=u_n$, for $n\in\enne$. We immediately see that
$\by(0)=y_0$ and that \eqref{weak1} holds true for $(\bm,\by)$. Also, the Lipschitz continuity of $f_2'$
and \eqref{p2bis} imply that $\,f_2'(\yn)\to f_2'(\by)\,$  strongly in \,$C^0([0,T];H)$.
Now, recall that $B^\sigma \yn\to B^\sigma\by$ weakly in $L^2(0,T;H)$, by virtue of \eqref{p2}. We thus have, by
\pier{lower} semicontinuity,
\begin{align}
&\int_0^T\left(B^\sigma\by(t),B^\sigma(\by(t)-v(t))\right)\, dt
  \non
  \\
& \leq \liminf_{n\to\infty} \ioT \bigl( B^\sigma \yn(t) , B^\sigma \yn(t) \bigr) \, dt
  - \lim_{n\to\infty} \ioT \bigl( B^\sigma \yn(t) , B^\sigma v(t) \bigr) \, dt
  \non
  \\
& = \liminf_{n\to\infty} \ioT \bigl( B^\sigma \yn(t) , B^\sigma (\yn(t) - v(t)) \bigr) \, dt
  \non
\end{align}
for every $v\in\L2{\VB\sigma}$. 
In conclusion, \gianni{owing to~\eqref{limh} as well}, we have that
\begin{align}
  & \intQ \indi(\by)
  + \ioT \bigl( B^\sigma \by(t) , B^\sigma (\by(t) - v(t)) \bigr) \, dt
  = \ioT \bigl( B^\sigma \by(t) , B^\sigma (\by(t) - v(t)) \bigr) \, dt
  \non
  \\
  & \leq \,\liminf_{n\to\infty} \intQ h^{\alpha_n}(\yn)
  + \liminf_{n\to\infty} \ioT \bigl( B^\sigma \yn(t) , B^\sigma (\yn(t) - v(t)) \bigr) \, dt
  \non
  \\
   & \leq \,\liminf_{n\to\infty} \Bigl(
    \intQ h^{\alpha_n}(\yn)
    + \ioT \bigl( B^\sigma \yn(t) , B^\sigma (\yn(t) - v(t)) \bigr) \, dt
  \Bigr)
  \non
  \\
  & \leq \,\lim_{n\to\infty} \,\Bigl( \ioT \bigl(
- \tau \partial_t y^{\alpha_n}(t) - f_2'(\yn(t)) + u(t) + \mn(t) , \yn(t) - v(t)
  \bigr)  \, dt
  + \intQ h^{\alpha_n}(v)\Bigr)
  \non
  \\
  & = \,\ioT \bigl(-\tau \partial_t \by(t) 
     -f_2'(\by(t)) + u(t) {\,+\,} \bm(t) , \by(t)-v(t) \bigr) \, dt
    + \intQ \indi(v),
  \non
\end{align}
for all $v\in L^2(0,T;V_B^\sigma)$. Thus the time-integrated version of \eqref{weak2}, with time-dependent
test functions, holds true. Since this version is equivalent to \eqref{weak2}, 
we see that $(\bm,\by)$ is indeed a solution in the sense of \cite[Thm.~2.6]{CGS18} to \eqref{weak1}--\eqref{weak3} for $f_1=\indi$. The assertion is thus proved. 
\Edim

\Brem
According to \cite[Thm.~2.6]{CGS18}, the second solution component $y$ and the expression $A^r\mu$ are uniquely determined.
This entails that $\by=\S_0(u)$ and that the convergence properties \pier{\eqref{p2}, \eqref{p3}} and \eqref{p5} are valid for the 
entire sequence $\{\alpha_n\}$ and not only for a subsequence. In addition, we can infer from \eqref{p1} that
$A^r\mn \to A^r\bm$ weakly in $L^2(0,T;H)$ as $n\to\infty$. If \,$\lambda_1>0$\pier{, then even $\mn$ converges to}  $\bm$ weakly in $L^2(0,T;V_A^r)$.  
\Erem
In the following theorem, we prove a quantitative estimate that yields information on the order of convergence as
$\alpha\searrow0$ in a very special (but important) situation. To this end, we need further assumptions that will
also be needed in the derivation of first-order necessary optimality conditions in Section~5.  

\Bthm
Suppose that in addition to \pier{{\bf (A1)}--{\bf (A6)}} the following assumptions are fulfilled:

\noindent
{\bf \gianni{(A8)}} \,\,\,The condition \eqref{adco} is satisfied.

\noindent
{\bf \gianni{(A9)}} \,\,\,$B=-\Delta$ with the domain $D(B)=\{v\in\Hdue: \dn v=0 \,\mbox{ on }\,\Gamma\}$, \,$\sigma=\frac 12$,
and \hspace*{12mm} $\,V_A^{2r}\subset L^\infty(\Omega)$.

\noindent
Moreover, assume that \,$u^{\alpha_1}, u^{\alpha_2}\in{\cal X}$ are given,
where \,$0<\alpha_1<\alpha_2<1$, and that $(\mu^{\alpha_i},y^{\alpha_i})$ are solutions to \eqref{weak1}--\eqref{weak3} 
for $f_1=h^{\alpha_i}$ and $u=u_i$ in the sense of {\rm \cite[Thm.~2.6]{CGS18}}, for
$i=1,2$.
Then there is a constant $K_2>0$, which
depends only on the data of the problem, such that it holds, for all $t\in (0,T]$, 
\begin{align}
\label{esti1}
&\|y^{\alpha_1}-y^{\alpha_2}\|_{C^0([0,t];\Ldue)\cap L^2(0,t;\Huno)}\,+\,\Bigl\|\int_0^{\bullet} A^r(\mu^{\alpha_1}-\mu^{\alpha_2})(s)ds\Bigr
\|_{C^0([0,t];\Ldue)}\non\\
&\le\,K_2\,\Bigl(|\alpha_1-\alpha_2|^{1/2}\,+\,\|u^{\alpha_1}-u^{\alpha_2}\|_{L^2(0,t;H)}\Bigr)\,.
\end{align}
\Ethm
\Bdim
We first observe that in \cite[Example~1]{CGS19} it has been shown that a uniform separation property is valid 
for the solutions to \eqref{weak1}--\eqref{weak3} with $f_1=h^\alpha$ under the assumptions {\bf (A1)--\gianni{(A9)}}; 
that is, there are constants $r_*,r^*\in (-1,1)$ \gianni{(depending on~$\alpha$)} such that
\begin{equation}
\label{separa}
  \gianni{r_* \le y^\alpha \le r^* \quad\mbox{a.e. in \,$Q$}}.
\end{equation}
Moreover, we have $V_{-\Delta}^{1/2}=\Huno$, and thus we can infer from \cite[\pier{Remarks}~3.4,~3.5~and~3.6]{CGS19} that
for any $\alpha>0$ the solution $(\mual,\yal)$ to \eqref{weak1}--\eqref{weak3} in the sense of \cite[Thm.~2.6]{CGS18} for
$f_1=h^\alpha$ is in fact uniquely determined and satisfies the variational equality (which in this special case
turns out to be equivalent to \eqref{weak2})
\begin{align}
\label{weak2new}
&(\tau\dt y^\alpha(t),v) + (\nabla \yal(t),\nabla v)+( (h^\alpha)'(\yal(t)),v)+(f_2'(\yal(t)),v)
=(\mual(t)+u(t),v)\nonumber\\
&\mbox{for a.e. $t\in (0,T)$ and every \,$v\in \Huno$.}
\end{align}   

Now\pier{,}  let $u:=u^{\alpha_1}-u^{\alpha_2}$, $\mu:=\mu^{\alpha_1}-\mu^{\alpha_2}$, and $y:=y^{\alpha_1}-y^{\alpha_2}$. 
Then, taking the difference 
in \eqref{weak1} for the two different cases $\alpha=\alpha_1$, $\alpha=\alpha_2$, and integrating the resulting
equality over $[0,t]$ with respect to time, where $t\in (0,T]$, we obtain the identity
\begin{equation*}
\langle y(t),v\rangle_{A,r}\,+\,\Bigl(A^r\mbox{$\int_0^t\mu(s)ds$},A^r v\Bigr)\,=\,0
\quad\mbox{for all $t\in (0,T]$ and $v\in V_A^r$}.
\end{equation*}
Testing this identity by $v=\mu(t)$, and noting that $\langle y(t),\mu(t)\rangle_{A,r}=(y(t),\mu(t))$ for almost
every $t\in(0,T)$, we thus obtain that
\begin{align}
\label{diff1}
&\int_0^t\!\!\int_\Omega y\mu\,=\,- \int_0^t\Bigl(A^r\mu(s), \int_0^s A^r\mu(\rho)d\rho\Bigr)ds
\,=\,-\frac 12\,\Bigl\|\int_0^t A^r\mu(s)ds\Bigr\|^2\,.
\end{align}
Next, we insert $v=\gianni{-y}$ in the variational equality \eqref{weak2new} for 
$\alpha=\alpha_2$, and $v=\gianni y$ in \eqref{weak2new} for $\alpha=\alpha_1$. Summation of the resulting
\pier{identities} then yields the equality
\begin{align}
\label{diff2}
&\frac\tau 2\,\|y(t)\|^2\,+\,\int_0^t\|\nabla y(s)\|^2\,ds\,+\int_0^t\!\!\int_\Omega \varphi(\alpha_1)(h'(y^{\alpha_1})-h'(y^{\alpha_2}))
(y^{\alpha_1}-y^{\alpha_2})\nonumber\\
&\gianni =\,-\int_0^t\!\!\int_\Omega(\varphi(\alpha_1)-\varphi(\alpha_2))h'(y^{\alpha_2})(y^{\alpha_1}-y^{\alpha_2})
\,-\int_0^t\!\!\int_\Omega(f_2'(y^{\alpha_1})-f_2'(y^{\alpha_2}))\,y\nonumber\\
&\,+\,\int_0^t\!\!\int_\Omega y\mu\,+\int_0^t\!\!\iO u\,y. 
\end{align}
Owing to the monotonicity of $\,h'$, the third summand on the left-hand side of \eqref{diff2} is nonnegative.
Moreover,   
$h'(y^{\alpha_2})(y^{\alpha_1}-y^{\alpha_2})\,\le\,h(y^{\alpha_1})-h(y^{\alpha_2})$ almost
everywhere in $Q$, since $\,h\in C^1(-1,1)\,$ is convex, and \,$\phi(\alpha_1)<\phi(\alpha_2)$. So
the first summand on the \rhs\ of \eqref{diff2}, which we denote by $I$, satisfies
\begin{align}
I\,&\le\,(\phi(\alpha_2)-\phi(\alpha_1))\int_0^t\!\!\iO(|h(y^{\alpha_1})|\,+\,|h(y^{\alpha_2})|)\,\le\,C_1\,
(\alpha_2-\alpha_1)\,,
\end{align}
with $\,C_1:=4\ln(2)\,|\Omega|\,T\,\|\phi'\|_{C^0([0,1])}\,$, where $|\Omega|$ denotes the volume of $\Omega$. 
Therefore, adding \eqref{diff1}
and \eqref{diff2}, and using the Lipschitz continuity of $f_2'$, we obtain from Young's inequality an 
estimate of the form
\begin{align}
&\frac\tau 2\,\|y(t)\|^2\,+\,\int_0^t\|\nabla  y(s)\|^2\,ds\,+\,
\frac 12\,\Bigl\|\int_0^t A^r\mu(s)ds\Bigr\|^2\non\\
&\le\,C_1\,|\alpha_1-\alpha_2|\,+\,(L+1)\int_0^t\!\!\iO|y|^2\,+\,\frac 14\int_0^t\!\!\iO|u|^2\,,
\end{align}
and \eqref{esti1} follows from Gronwall's lemma.
\Edim

\Bcor
Suppose that {\bf (A1)}--{\bf \gianni{(A9)}} are fulfilled and that $\alpha\in (0,1]$. 
Moreover, let $y=\S_0(u)$ and $\yal=\S_{\alpha}(u)$. Then
\begin{align}
\label{esti2}
&\|y^{\alpha}-y\|_{C^0([0,t];\Ldue)\cap L^2(0,t;\Huno)}\,\le\,K_2\,|\alpha|^{1/2}\,.
\end{align}
\Ecor

\Bdim
We apply \eqref{esti1} with $\alpha_1=\alpha$, $\alpha_2=\alpha_n$, where $\alpha_n\searrow0$, and
$u^{\alpha_1}=u^{\alpha_2}=u$, which with $\yn=\S_{\alpha_n}(u)$ yields the estimate
\begin{align*}
&\|y^{\alpha}-\yn\|_{C^0([0,t];\Ldue)\cap L^2(0,t;\Huno)}\,+\,\Bigl\|\int_0^{\bullet} A^r(\mu^{\alpha}-\mn
)(s)ds\Bigr
\|_{C^0([0,t];\Ldue)}\non\\
&\le\,K_2\,|\alpha-\alpha_n|^{1/2}\,.
\end{align*}
The assertion now follows from \eqref{cony} in Theorem 3.2  by taking the
limit as $n\to\infty$, invoking Remark 3.4 
 and the semicontinuity of norms.
\Edim

\section{Existence and approximation of optimal controls}

\setcounter{equation}{0}
Beginning with this section, we investigate the optimal control problem \CP0 of minimizing the cost 
functional \eqref{cost} over the admissible set $\uad$ subject to state system \eqref{weak1}--\eqref{weak3} 
where $f_1=\indi$. In comparison with \CP0, 
we consider for $\alpha>0$ the following control problem:

\vspace{1mm}\noindent
\CPal \,\,\,Minimize $\,{\cal J}(y,u)\,$
for $\,u\in\uad$, subject to the condition that $y=\S_\alpha(u)$ for some solution $(\mu,y)$ to the state system \eqref{weak1}--\eqref{weak3} with $f_1=h^\alpha$\pier{,} in the sense of \cite[Thm.~2.6]{CGS18}.

\vspace{1mm}\noindent
We expect that the minimizers  of \CPal\ are for $\alpha\searrow0$ related to minimizers of \CP0.
Prior to giving an affirmative answer to this conjecture, we first show an existence result for \CPal.

\Bprop
Suppose that \pier{{\bf (A1)}--{\bf (A6)}} are satisfied. Then \mbox{{\rm(}${\cal CP}_{\alpha}${\rm)}} has for
every $\alpha>0$ a solution.
\Eprop

\Bdim
Let $\alpha>0$ be fixed, and assume that a  minimizing sequence $\,\{(y_n,u_n)\}$ for \CPal\ is given, 
where $y_n=\S_\alpha(u_n)$ for some solution
pair $(\mu_n,y_n)$ to the state system with $u=u_n\in\uad$ and $f_1=h^\alpha$, for $n\in\enne$. Then it holds for every $n\in\enne$ that
\begin{align}
\label{exi1}
&\langle \dt y_n(t),v\rangle_{A,r}+(A^r\mu_n(t),A^r v)=0\quad \mbox{for a.e. \,$t\in (0,T)$\, and every $\,v\in V_A^r$,}
\\[0.5mm]
\label{exi2}
&(\tau\dt y_n(t),y_n(t)-v)+(B^\sigma y_n(t),B^\sigma(y_n(t)-v))+h^\alpha(y_n(t))\non\\
&\quad \le\,(\mu_n(t)+u_n(t)-f_2'(y_n(t)),y_n(t)-v) +h^\alpha(v) \non\\
& \mbox{for a.e. \,$t\in (0,T)$\, and every $\,v\in V_B^\sigma$,}
\\[0.5mm]
\label{exi3}
&y_n(0)=y_0.
\end{align}
Taking \gianni{estimate \eqref{albound1}} into account, we may without loss of generality assume that there
are $\,\bu\in\uad\,$ and $(\bm,\by)$ such that
\begin{align}
\label{exiu}
u_n\to\bu&\quad\mbox{weakly-star in }\,{\cal X},\\
\label{eximu}
\mu_n\to\bm&\quad\mbox{weakly in }\,L^2(0,T;V_A^r),\\
\label{exiy}
y_n\to\by&\quad\mbox{weakly-star in }\,H^1(0,T;V_A^{-r})\cap L^\infty(0,T;V_B^\sigma),\nonumber\\
&\quad\mbox{strongly in }\,C^0([0,T];H), \,\,\mbox{ and pointwise a.e. in }\,\,Q.
\end{align} 
Then also $\,f_2'(y_n)\to f_2'(\by)\,$ strongly in $C^0([0,T];H)$. Moreover, it holds 
\begin{equation}\label{elvis}
\int_0^T \!\!\iO h^\alpha(y_n)\,\le\,K_1\quad\mbox{for every $n\in\enne$}.
\end{equation}
Therefore, we have $\,y_n\in [-1,1]\,$ almost everywhere in $Q$, and since $h^\alpha$ is continuous
in $[-1,1]$, it follows that $\,h^\alpha(y_n)\to h^\alpha(\by)$\, pointwise almost everywhere in $Q$.
Lebesgue's dominated convergence theorem then yields that
$$\int_0^T\!\!\iO h^\alpha(y_n)\to \int_0^T\!\!\iO h^\alpha(\by)\,.$$
In addition, by lower semicontinuity, \pier{we have that}
$$\int_0^T (B^\sigma \by(t)),B^\sigma \by(t))\,dt\,\le\,\liminf_{n\to\infty}\int_0^T
(B^\sigma y_n(t),B^\sigma y_n(t))\,dt\,.$$
Combining the convergence results shown above, we obtain by passage to the limit as $n\to\infty$ that 
\begin{align}
\label{exi4}
&\int_0^T\langle \dt \by(t),v(t)\rangle_{A,r}\,dt\,+\int_0^T(A^r\bm(t),A^r v(t))\,dt\,=\,0\quad\forall\,
v\in L^2(0,T;V_A^r),
\\[0.5mm]
\label{exi5}
&\int_0^T(\tau\dt \by(t),\by(t)-v(t))\,dt\,+\int_0^T(B^\sigma \by(t),B^\sigma(\by(t)-v(t)))\,dt
\,+\int_0^T\!\!\iO h^\alpha(\by)\non\\
&\le\int_0^T(\bm(t)+\bu(t)-f_2'(\by(t)),\by(t)-v(t))\,dt\,+\int_0^T\!\!\iO h^\alpha(v)\non\\[0.5mm] 
& \mbox{for every $\,v\in L^2(0,T;V_B^\sigma)$,}
\\[0.5mm]
\label{exi6}
&\by(0)=y_0.
\end{align}
Apparently, \eqref{exi4}--\eqref{exi5} is just the time-integrated version of 
\eqref{weak1}--\eqref{weak2} for $u=\bu$ and $f_1=h^\alpha$, 
written with time-dependent test functions, which is equivalent to 
\eqref{weak1}--\eqref{weak2}.
Hence, $(\bm,\by)$ solves \eqref{weak1}--\eqref{weak2} for $u=\bu$ and $f_1=h^\alpha$ in the sense of
\cite[Thm.~2.6]{CGS18}. In particular, we have
$\by=\S_\alpha(\bu)$. But this means that $\,(\by,\bu)\,$ is  admissible
for \CPal. From the semicontinuity properties
of the cost functional \eqref{cost} it then follows that $\,(\by,\bu)\,$ 
is an optimal pair, which concludes the proof of the assertion.
\Edim

\Bprop
Let \pier{{\bf (A1)}--{\bf (A6)}} be fulfilled, and suppose that sequences $\,\{\alpha_n\}\subset (0,1]\,$ and
$\,\{u_n\}\subset\uad\,$ are given such that $\,\alpha_n\searrow0\,$ and $\,u_n\to u\,$ weakly-star in ${\cal X}$
for some $\,u\in\uad$. Then, with the solution operators defined in \eqref{defcs0} and \eqref{defcsal},
\begin{align}\label{cesareuno}
&{\mathcal{J}}(\S_0(u),u)\,\le\,\liminf_{n\to\infty}\,{\mathcal{J}}(\S_{\alpha_n}(\un),\un),\\[0.5mm]
\label{cesaredue}
&{\mathcal{J}}(\S_0(v),v)\,=\,\lim_{n\to\infty}\,{\mathcal{J}}(\S_{\alpha_n}(v),v) \quad\forall\,v\in\uad.
\end{align}
\Eprop

\Bdim
Under the given assumptions, we may apply \eqref{cony}  in Theorem 3.3 and Remark 3.4 to infer
that $\gianni{\S_{\alpha_n}(\un)\to\S_0(u)}$ strongly in $C^0([0,T];H)$.  
The validity of \eqref{cesareuno} is then a direct consequence of the weak and weak-star 
sequential semicontinuity properties of the cost
functional~\gianni{\eqref{cost}}. 
Now suppose that $v\in\uad$ is arbitrarily chosen, and put $\yn:=\S_{\alpha_n}(v)$ for all $n\in\enne$.
Then, again by Theorem~3.3 and Remark~3.4, $\,\yn\to \S_0(v)\,$ strongly in $C^0([0,T];H)$. \pier{Next,} 
observe that the first two summands of the cost functional are 
obviously continuous with respect to the strong topology of $C^0([0,T];H)$, which then shows the validity of
\eqref{cesaredue}. 
\Edim

We are now in a position to prove the existence of minimizers for the problem \CP0. We have the following result.

\Bcor
Under the assumptions of Proposition~4.2, the optimal control problem {\rm ($\mathcal{CP}_0$)} has at least
one solution.
\Ecor

\Bdim
Pick an arbitrary sequence $\{\alpha_n\}\subset (0,1]$ such that $\alpha_n\searrow0$ as $n\to\infty$.
By virtue of Proposition 4.1, the optimal control problem (${\mathcal{CP}}_{\alpha_n}$) has for every $n\in\enne$ a solution $(\yn,\un)$, where $\yn=\S_{\alpha_n}(\un)$ for a solution $(\mn,\yn)$ to the corresponding state system.
Since $\uad$ is bounded in ${\cal X}$, we may without loss of generality assume that $\un\to u$ weakly-star in ${\cal X}$
for some $u\in\uad$.  At this point, we apply Theorem 3.3 to the present situation. We then infer that the convergence
results \eqref{conmu} and \eqref{cony} hold true for some subsequence $\,\{\alpha_{n_k}\}\,$ 
with a pair $(\mu,y)$ satisfying $y=\S_0(u)$. Invoking the optimality of
$\,(\yn,\un)\,$ for (${\mathcal{CP}}_{\alpha_n}$), we then find for every $\,v\in\uad\,$ the chain of (in)equalities
\begin{align}
&{\cal J}(y,u)\,=\,{\cal J}(\S_0(u),u)\,\le\,\liminf_{k\to\infty}\,{\cal J}(\S_{\alpha_{n_k}}(u^{\alpha_{n_k}}),
u^{\alpha_{n_k}})\non\\
&\le\,\liminf_{k\to\infty}\,{\cal J}(\S_{\alpha_{n_k}}(v),v)\,\le\,{\cal J}(\S_0(v),v),
\end{align}
which yields that $\,(y,u)\,$ is an optimal pair for \CP0. The assertion is thus proved.
\Edim
Theorem 3.3 and the proof of Corollary 4.3 indicate that optimal controls of \CPal\ are ``close'' to optimal
controls of \CP0. However, they do not yield any information on whether every optimal control
of \CP0\ can be approximated in this way. In fact, such a global result cannot be expected to hold true. 
However, a local answer can be given. For this purpose, we employ a trick introduced in \cite{Barbu}. To this end, 
let $\bu\in\uad$ be an optimal control for \CP0\ with the associated state $(\bm,\by)$ where $\by=\S_0(\bu)$. We associate
with this optimal control the {\em adapted cost functional}
\begin{equation}
\label{adcost}
\widetilde{\cal J}(y,u):={\cal J}(y,u)\,+\,\frac 12\,\|u-\bu\|^2_{L^2(Q)}
\end{equation}
and a corresponding \emph{adapted optimal control problem} for $\alpha>0$, namely:

\vspace{2mm}\noindent
($\widetilde{\mathcal{CP}}_{\alpha}$)\quad Minimize $\,\, \widetilde {\cal J}(y,u)\,\,$
for $\,u\in\uad$, subject to the condition that $y=\S_\alpha(u)$ for some solution $(\mu,y)$ to the state system \eqref{weak1}--\eqref{weak3} with $f_1=h^\alpha$ in the sense of \cite[Thm.~2.6]{CGS18}.

\vspace{2mm}
With essentially the same proof as that of Proposition 4.1 (which needs no repetition here), we can show the following 
result.

\Blem
Suppose that the assumptions \pier{{\bf (A1)}--{\bf {(A6)}}} are fulfilled. Then the optimal control problem 
$(\widetilde{\cal CP}_\alpha)$ has for every $\alpha>0$ at least one solution.
\Elem
  
\vspace{1mm}
We are now in the position to give a partial answer to the question raised above. We have the following result.

\Bthm
Let the assumptions \pier{{\bf (A1)}--{\bf {(A6)}}} be fulfilled, suppose that 
$\bar u\in \uad$ is an arbitrary optimal control of {\rm $({\mathcal{CP}}_{0})$} with associated state  
$(\bar\mu,\by)$ where $\by=\S_0(\bu)$, and let $\,\{\alpha_n\}\subset (0,1]\,$ be
any sequence such that $\,\alpha_n\searrow 0\,$ as $\,n\to\infty$. Then there exist a 
subsequence $\{\alpha_{n_k}\}$ of $\{\alpha_n\}$, and, for every $k\in\enne$, an optimal control
 $\,u^{\alpha_{n_k}}\in \uad\,$ of the adapted problem {\rm $(\widetilde{\mathcal{CP}}_{\alpha_{n_k}})$}
 with associated state $(\mu^{\alpha_{n_k}},y^{\alpha_{n_k}})$, where $\,y^{\alpha_{n_k}}=\S_{\alpha_{n_k}}(u^{\alpha_{n_k}})$,
 such that, as $k\to\infty$,
\begin{align}
\label{tr3.4}
&u^{\alpha_{n_k}}\to \bar u\quad\mbox{strongly in }\,L^2(Q),
\end{align}
and such that the property \eqref{cony} is satisfied with  $\,y\,$  replaced
by  $\,\by\,$. Moreover, we have 
\begin{align}
\label{tr3.5}
&\lim_{k\to\infty}\,\widetilde{{\cal J}}(y^{\alpha_{n_k}},u^{\alpha_{n_k}})\,=\,  {\cal J}(\bar y,\bar u)\,.
\end{align}
\Ethm

\Bdim
Let $\alpha_n \searrow 0$ as $n\to\infty$. For any $ n\in\enne$, we pick an optimal control 
$u^{\alpha_n} \in \uad\,$ for the adapted problem $(\widetilde{\cal CP}_{\alpha_n})$ and denote by 
$(\mu^{\alpha_n},\yn)$, where $\yn=\S_{\alpha_n}(\un)$, an  associated solution to  
(\ref{weak1})--(\ref{weak3}) with $f_1=h^{\alpha_n}$ and $u=\un$. 
By the boundedness of $\uad$ in $\calX$, there is some subsequence $\{\alpha_{n_k}\}$ of $\{\alpha_n\}$ such that
\begin{equation}
\label{ugam}
u^{\alpha_{n_k}}\to u\quad\mbox{weakly-star in }\,{\cal X}
\quad\mbox{as }\,k\to\infty,
\end{equation}
with some $u\in\uad$. 
Thanks to Theorem~3.3, the convergence \gianni{properties \accorpa{conmu}{cony} hold} true 
with some pair $(\mu,y)$ satisfying $y=\S_0(u)$. In particular,  the pair $(y,u)$
is admissible for (${\cal CP}_0$). 

We now aim to prove that $u=\bar u$. Once this is shown, it follows from the uniqueness of the second
solution component to the state system
\eqref{weak1}--\eqref{weak3} that also $y=\bar y$, 
which implies that \gianni{\eqref{cony}} holds true with  $\,y\,$ replaced by  $\,\by$.

Now observe that, owing to the weak sequential lower semicontinuity of 
$\widetilde {\cal J}$, 
and in view of the optimality property of $(\by,\bar u)$ for problem $({\cal CP}_0)$,
\begin{align}
\label{tr3.6}
&\liminf_{k\to\infty}\, \widetilde{\cal J}(y^{\alpha_{n_k}},
u^{\alpha_{n_k}})
\ge \,{\cal J}(y,u)\,+\,\frac{1}{2}\,
\|u-\bar{u}\|^2_{L^2(Q)}\nonumber\\[1mm]
&\geq \, {\cal J}(\by,\bar u)\,+\,\frac{1}{2}\,\|u-\bar{u}\|^2_{L^2(Q)}\,.
\end{align}
On the other hand, the optimality property of  $\,(y^{\alpha_{n_k}},u^{\alpha_{n_k}})
\,$ for problem $(\widetilde {\cal CP}_{\alpha_{n_k}})$ yields that
for any $k\in\enne$ we have
\begin{equation}
\label{tr3.7}
\widetilde {\cal J}(y^{\alpha_{n_k}},u^{\alpha_{n_k}})\, =\,
\widetilde {\cal J}({\cal S}_{\alpha_{n_k}}(u^{\alpha_{n_k}}),
u^{\alpha_{n_k}})\,\le\,\widetilde {\cal J}({\cal S}_{\alpha_{n_k}}
(\bar u),\bar u)\,,
\end{equation}
whence, taking the limit superior as $k\to\infty$ on both sides and invoking (\ref{cesaredue}) in
Proposition 4.2,
\begin{align}
\label{tr3.8}
&\limsup_{k\to\infty}\,\widetilde {\cal J}(y^{\alpha_{n_k}},
u^{\alpha_{n_k}})\,\le\,\widetilde {\cal J}({\calS}_0(\bar u),\bar u) 
\,=\,\widetilde {\cal J}(\bar y,\bar u)
\,=\,{\cal J}(\bar y,\bar u)\,.
\end{align}
Combining (\ref{tr3.6}) with (\ref{tr3.8}), we have thus shown that 
$\,\frac{1}{2}\,\|u-\bar{u}\|^2_{L^2(Q)}=0$\,,
so that $\,u=\bar u\,$  and thus also $\,y=\bar y$. 
Moreover, (\ref{tr3.6}) and (\ref{tr3.8}) also imply that
\begin{align}
\label{tr3.9}
&{\cal J}(\bar y,\bar u) \, =\,\widetilde{\cal J}(\bar y,\bar u)
\,=\,\liminf_{k\to\infty}\, \widetilde{\cal J}(y^{\alpha_{n_k}},
 u^{\alpha_{n_k}})\nonumber\\[1mm]
&\,=\,\limsup_{k\to\infty}\,\widetilde{\cal J}(y^{\alpha_{n_k}},
u^{\alpha_{n_k}}) \,
=\,\lim_{k\to\infty}\, \widetilde{\cal J}(y^{\alpha_{n_k}},
u^{\alpha_{n_k}})\,,
\end{align}                                     
which proves {(\ref{tr3.5})} and, at the same time, also (\ref{tr3.4}). This concludes the proof
of the assertion.
\Edim

\section{Adjoint system and first-order optimality conditions}

\setcounter{equation}{0}

In this section, we aim at deriving first-order necessary optimality conditions for the optimal
control problem \CP0 using the deep quench approximation. Throughout the section, we assume
that $\bu\in\uad$ is an optimal control of (${\mathcal{CP}}_0$) with associated state
$(\bm,\by)$\pier{, with $\by=\S_0(\bu)$.}  The derivation will be achieved by a passage to the 
limit as $\alpha\searrow0$ in the first-order optimality conditions for the adapted optimal control
problems ($\widetilde{\mathcal{CP}}_\alpha$) that can be derived as in \cite{CGS19} with only minor
and obvious changes.  This approach will not be possible in full generality. In fact, we have to assume
that, besides {\bf (A1)--\gianni{(A7)}}, the assumptions {\bf \gianni{(A8)}--\gianni{(A9)}} from Theorem 3.5 are fulfilled.

\Brem
Observe that {\bf \gianni{(A8)}} yields the validity of the stronger regularity properties \eqref{albound2} from Theorem~3.1. 
Also, {\bf \gianni{(A9)}} implies that the constant functions belong to $V_{-\Delta}^{1/2}=\Huno$, so that {\bf (A2)}
is automatically fulfilled. 
In addition, since $\Huno\cap L^\infty(\Omega)$ is dense in $\Huno$ and 
$\Huno$ is continuously embedded in $L^4(\Omega)$, the conditions \cite[{\bf (A8)}~and~{\bf (A9)}]{CGS19} are satisfied. 
\Erem
\Brem
The condition that $\,V_A^{2r}\subset L^\infty(\Omega)\,$ is, for instance, 
satisfied if $A=-\Delta$ with zero Dirichlet or Neumann
 boundary condition and $r> \frac 3 8$. Indeed, we have in this case that   $V_A^{2r}
\subset H^{4r}(\Omega)\subset L^\infty(\Omega)$, since $\,\,4r>\frac 32$. Likewise, if $A=\Delta^2$ with
domain $D(A)\subset H^4(\Omega)$, then $V_A^{2r}\subset L^\infty(\Omega)$ provided that $r>\frac 3{16}$.
In this sense, while the 
improvement obtained in the following results over previously known results for the 
classical case $A=B=-\Delta$, $r=\sigma=
\frac 12$, is not too large, the results are entirely new for other operators \pier{$A$}; in fact, to our best knowledge, 
they constitute the first ever first-order necessary
optimality conditions for Cahn--Hilliard type systems with fractional operators and nondifferentiable nonlinearities of 
double obstacle type.  
\Erem

As already mentioned in the proof of Theorem 3.5, it follows under the assumptions {\bf (A1)--\gianni{(A9)}} that also the
solution component $\,\mual\,$ of the solutions $(\mual,\yal)$  to \eqref{weak1}--\eqref{weak3} in the sense of \cite[Thm.~2.6]{CGS18} for $f_1=h^\alpha$ is uniquely determined, so that a corresponding solution operator
$$\widetilde\S_\alpha=(\widetilde \S_\alpha^1,\widetilde\S_\alpha^2):u\ni\uad\mapsto\widetilde\S_\alpha(u)=
(\widetilde \S_\alpha^1(u),\widetilde\S_\alpha^2(u)):=(\mual,\yal)$$
is well defined. Clearly, we have $\widetilde\S_\alpha^2=\S_\alpha$.  Moreover, $(\mual,\yal)$ satisfies the variational equality \eqref{weak2new}, which in this situation is equivalent to \eqref{weak2}. In addition, a 
uniform separation property is satisfied; indeed, \pier{thanks to {\bf (A6)}}, for every $\alpha>0$ and every bounded set ${\cal U}\subset {\cal X}$, 
there exist constants $r_*(\alpha),r^*(\alpha)\in (-1,1)$, which  depend only on ${\cal U}$,  
such that the following holds true: whenever $(\mual,\yal)=\widetilde\S_\alpha(u)$ 
for some $u\in {\cal U}$, then
\begin{equation}
\label{howgh}
r_*(\alpha) \le y^{\alpha} \le r^*(\alpha) \quad\mbox{a.e. in \,$Q$}, \quad r_*(\alpha)\le y_0\le r^*(\alpha)\quad
\mbox{a.e. in $\,\Omega$.}
\end{equation}
In particular, the condition \cite[{\bf (GB)}]{CGS19}, which was crucial for the analysis carried out in 
\cite{CGS19}, is fulfilled for the potentials $f_1=h^\alpha$, $\alpha>0$, and we may take advantage of the results derived there.
\Brem  Owing to the separation property \eqref{howgh}, there is, for every $\alpha>0$ and every
bounded ${\cal U}\subset {\cal X}$, some constant $K_\alpha>0$, which depends only on ${\cal U}$, such
that
\begin{equation}
\label{albound3}
\max_{0\le i\le 3}\,\|(h^\alpha)^{(i)}(\yal)\|_{L^\infty(Q)}\,\le\,K_\alpha \quad\mbox{whenever $\,\yal=\S_\alpha(u)\,$
for some $\,u\in {\cal U}$}.
\end{equation}
Now we have  $V_A^{2r}
\subset L^\infty(\Omega)$ and thus, by  \eqref{albound2},     $\mual\in L^\infty(Q)$. Since also 
$\dt\yal\in L^\infty(0,T;H)$, comparison in \eqref{weak2new} shows that then $\yal\in L^\infty(0,T;\Hdue)$, which means
that the state equations \eqref{weak1}, \eqref{weak2} for $f_1=h^\alpha$ are even satisfied in the strong sense, 
that is, we have
\begin{align}
\label{pw1}
&\dt\yal+A^{2r}\mual=0 \quad\mbox{a.e. in }\,Q,\\
\label{pw2}
&\tau\dt \yal-\Delta\yal+(h^\alpha)'(\yal)+f_2'(\yal)=\mual+u \quad\mbox{a.e. in }\,Q.
\end{align}
\Erem
At this point, we observe that that the state systems associated with (${\mathcal{CP}}_\alpha$) and 
($\widetilde{\mathcal{CP}}_\alpha$) are exactly the same. Hence, if $\bual\in\uad$ is an optimal control
of ($\widetilde{\mathcal{CP}}_\alpha$) with associated state $(\bmal,\byal)=\widetilde\S_\alpha(\bual)$ 
for some $\alpha>0$, then 
$(\bmal,\byal)$ satisfies the global bounds \eqref{albound1}, \eqref{albound2}, \eqref{albound3}, as well as the
separation property \eqref{howgh}, and the state equations hold true in the form \eqref{pw1}, \eqref{pw2}. Moreover,
introducing for $\alpha>0$ the abbreviating notation
\begin{equation}
\label{defgal}
g_1^\alpha:=\beta_1(\byal(T)-y_\Omega),\quad g_2^\alpha:=\beta_2(\byal-y_Q),\quad \psi_1^\alpha:=f_2''(\byal),
\quad\psi_2^\alpha:=\phi(\alpha)
h''(\byal),
\end{equation}
we observe that \eqref{albound1}, \eqref{albound2}, \eqref{interv}, and {\bf \gianni{(A5)}} imply the global bound
\begin{equation}
\label{galbound}
\|g_1^\alpha\|_{\Ldue}\,+\,\|g_2^\alpha\|_{L^2(Q)}\,+\,\|\psi_1^\alpha\|_{L^\infty(Q)}\,\le\,C_1\quad\forall\,\alpha \in (0,1],
\end{equation}
where, here and in the following, $C_i$, $i\in\enne$, denote positive constants that may depend on the data of the state
system but not on $\alpha\in (0,1]$. 
Observe that a corresponding bound for $\psi^\alpha_2$ cannot be expected: indeed, it may well happen that 
the separation constants $r_*(\alpha)$ and/or $r^*(\alpha)$ introduced in \eqref{howgh} approach $\pm 1$ as
$\alpha\searrow0$, so that $\,\psi_2^\alpha=\frac {2\phi(\alpha)}{1-(\byal)^2}\,$ may become unbounded as
$\alpha\searrow0$.

\vspace{3mm}
Next, we consider the adjoint system associated with the 
adapted optimal control problem ($\widetilde{\mathcal{CP}}_\alpha$). According to \cite[Sect.~5]{CGS19}, it has the  
following form:
\begin{align}
\label{adj1}
&(A^r\pal(t),A^r v)-(\qal(t),v)\,=\,0\quad\mbox{for a.e. $\,t\in (0,T)\,$ and every }\,v\in V_A^r,\\[0.5mm]
\label{adj2} 
&\langle -\dt (\pal+\tau\qal)(t),v\rangle +(\nabla \qal(t),\nabla v)+((\psi_1^\alpha(t)+\psi_2^\alpha(t))\,\qal(t),v)\non\\
&=\,(g_2^\alpha(t),v) \quad\mbox{for a.e. $\,t\in (0,T)\,$ and every }\,v\in 
\Huno,\\[0.5mm]
\label{adj3}
&(\pal+\tau\qal)(T)\,=\,g_1^\alpha \quad\mbox{in }\,\Omega.
\end{align} 
Here, for the sake of simplicity, we have denoted by $\langle\cdot,\cdot\rangle$ the dual pairing between
${\Huno}^*$ and $\Huno$. 

The system \eqref{adj1}--\eqref{adj3} is a special case of the type of systems that has been analyzed in
\cite[Sect.~5]{CGS19}. We briefly summarize some of the results established there 
(cf., \cite[Prop.~5.2,~Lem.~5.3,~Lem.~5.4, Rem.~5.7,~Thm.~5.8]{CGS19}), 
 where we have to distinguish the following cases:

\vspace{2mm}\noindent
\underline{{\sc Case 1:} $\lambda_1>0$.}\\[1mm]
In this case, the system \eqref{adj1}--\eqref{adj3}
 admits a unique solution $(\pal,\qal)$ satisfying
\begin{align}
\label{regp}
&\pal\in L^2(0,T;V_A^{2r}),\\
\label{req}
&\qal\in L^2(0,T;\Huno),\\
\label{regpq}
&\pal+\tau \qal\in H^1(0,T;{\Huno}^*).
\end{align}
Notice that \eqref{regpq} implies that $p+\tau q\in C^0([0,T]; {\Huno}^*)$, so that the endpoint condition 
\eqref{adj3} is meaningful.
Now observe that the operator $A^{2r}\in {\cal L}(V_A^{2r},H)$ is for $\lambda_1>0$ a topological isomorphism, 
and  with $\,A^{-2r}:=(A^{2r})^{-1}:
H\to V_A^{2r}$ the variational equation \eqref{adj1} takes the simple form $\,\pal=A^{-2r}\qal$. Inserting this in
\eqref{adj2} and \eqref{adj3}, we obtain that
\begin{align}
\label{lg1}
&\langle -\dt \bigl((A^{-2r}+\tau I)\qal\bigr)(t),v\rangle +\iO\nabla \qal(t)\cdot\nabla v\,+(\psi_1^\alpha(t)\,\qal(t),v)\non
\\[0.5mm]
&\quad+(\psi_2^\alpha(t)\qal(t),v)\,=\,(g_2^\alpha(t),v) \quad\mbox{for a.e. $\,t\in (0,T)\,$ and all }\,v\in 
\Huno,\\[1mm]
\label{lg2}
&(A^{-2r}+\tau I)\qal(T)\,\,=\,g_1^\alpha\qquad\mbox{a.e. in }\,\Omega,
\end{align} 
where $\,I\,$ denotes the identity operator in $\,H$. 
Moreover, since also the linear operator $\,A^{-2r}+\tau I \in {\cal L}(H,H)$ is obviously
a topological isomorphism,  \eqref{lg2} can be equivalently written as
\begin{equation}
\label{lg2neu}
\qal(T)\,=\,(A^{-2r}+\tau I)^{-1} g_1^\alpha,
\end{equation} 
which gives $\qal(T)$ a proper meaning as well. 

We now derive an estimate for the adjoint variables that is uniform in $\alpha>0$.
Testing \eqref{lg1} by $\qal(t)$ and integrating with respect to time over $[t,T]$, where $t\in[0,T)$, we then
conclude the equation
\begin{align}
\label{lg3}
&\int_t^T\!\langle -\dt \bigl((A^{-2r}+\tau I)\qal\bigr)(\rho),\qal(\rho)\rangle\,d\rho\, 
+\int_t^T\!\|\nabla \qal(\rho)\|^2\,d\rho
\, +\int_t^T\!\!\iO\!\psi_2^\alpha \,|\qal|^2\non
\\[0.5mm]
&=\,\int_t^T\!\!\iO\bigl(-\psi_1^\alpha \qal\,+\,g_2^\alpha\bigr)\,\qal\,,
\end{align}
where the last term on the \lhs\ is nonnegative and, owing to \eqref{galbound}, the \rhs\ is bounded by an expression of the form
\begin{equation}
\label{Paul}
C_2\,+\,C_3\int_t^T\!\!\iO|\qal|^2\,.
\end{equation}
Now observe that, by definition \eqref{defAr}, and since $\lambda_1>0$, it holds for every $v\in H$ that
\begin{equation}
\label{lg4}
(A^{-2r}+\tau I)^{1/2}v\,=\,\sum_{j=1}^\infty \left(\lambda_j^{-2r}+\tau\right)^{1/2}(v,e_j)e_j,
\end{equation}
and we have the estimates
\begin{align}
  \label{lg5}
  &\left\|(A^{-2r}+\tau I)^{1/2}v\right\|^2\,=\,\sum_{j=1}^\infty\left(\lambda_j^{-2r}+\tau\right)|(v,e_j)|^2
  \,\ge\,\tau\,\|v\|^2,
  \\
  \label{lg6}
  &\left\|(A^{-2r}+\tau I)^{1/2}v\right\|^2\,\gianni{{=\,\sum_{j=1}^\infty\left(\lambda_j^{-2r}+\tau\right)|(v,e_j)|^2}}
  \,\le\,\gianni{(\lambda_1^{-2r}+\tau)}\,\|v\|^2,
  \\
  \label{lg7}
  &\left\|(A^{-2r}+\tau I)^{-1}v\right\|^2\,=\,\sum_{j=1}^\infty(\lambda_j^{-2r}+\tau)^{-2}|(v,e_j)|^2
  \,\le\,\tau^{-2}\,\|v\|^2.
\end{align} 
Moreover, it is easily verified that
\begin{align}
\label{lg8}
-\,\langle \dt(A^{-2r}+\tau I)\qal(t),\qal(t)\rangle\,=\,-\,\frac 12\,\frac d{dt}\left\|(A^{-2r}+\tau I)^{1/2}
\qal(t)\right\|^2.
\end{align}
Therefore, by virtue of \eqref{lg2neu}, the first term on the \lhs\ of \eqref{lg3} is equal to
\begin{equation}
\label{lg9}
\frac 12\left\|(A^{-2r}+\tau I)^{1/2}\qal(t)\right\|^2\,-\,\frac 12\left\|(A^{-2r}+\tau I)^{1/2}
(A^{-2r}+\tau I)^{-1} g_1^\alpha\right\|^2,
\end{equation}
which, by \eqref{galbound} and \eqref{lg5}--\eqref{lg7}, is bounded from below by $\,\frac \tau 2 \,\|\qal(t)\|^2-C_4$,
with some global constant $C_4>0$. At this point, we invoke Gronwall's lemma, taken backward in time,
as well as the fact that $p=A^{-2r}q$, to conclude that
\begin{equation}
\label{lgq}
\|\pal\|_{L^\infty(0,T;V_A^{2r})}\,+\,\|\qal\|_{L^\infty(0,T;H)\cap L^2(0,T;\Huno)}\,\le\,C_5\quad\forall\,
\juerg{\alpha\in (0,1]}.
\end{equation}

\vspace{2mm}\noindent
\underline{{\sc Case 2:} $\lambda_1=0.$}\\[1mm]
This case is considerably more difficult to handle. To motivate this, we 
denote by $\,\mathbf{1}\,$ both the functions that are identically equal to  \,$1$\, 
in either $\,\Omega\,$ or $\,Q$.
Then, by {\bf (A2)}(ii), $\,A^r{\mathbf{1}}=0$, and insertion
of $v={\mathbf{1}}$ in \eqref{adj1} yields that
\Beq
\label{meanqal}
\mbox{mean}\,(\qal(t))=0 \quad\mbox{for a.e. $\,t\in (0,T)$.}
\Eeq

At this point, and also for later use, we recall an integration-by-parts formula that was proved in
\cite[Lem.~4.5]{CGSSIAM18}: if $({\cal V},{\cal H},{\cal V}^*)$ is a Hilbert triple and
\begin{equation}
\label{leibniz1}
w\in H^1(0,T;{\cal H})\cap L^2(0,T;{\cal V})\quad\mbox{and}\quad z\in H^1(0,T;{\cal V}^*)\cap L^2(0,T;{\cal H}),
\end{equation}
then the function $\,t\mapsto (w(t),z(t))_{\cal H}\,$ is absolutely continuous, and for every $t_1,t_2\in [0,T]$ it holds 
the formula
\begin{equation}
\label{leibniz2}
\int_{t_1}^{t_2}\bigl[(\dt w(t),z(t))_{\cal H}\,+\,\langle \dt z(t),w(t)\rangle_{{\cal V}}\bigr]\,dt
\,=\,(w(t_2),z(t_2))_{{\cal H}}-(w(t_1),z(t_1))_{{\cal H}},
\end{equation}
where $\,(\cdot,\cdot)_{\cal H}\,$ denotes the inner product in ${\cal H}$.

We now insert $v=\mathbf{1}$ in \eqref{adj2} and integrate the resulting identity with respect to time
over $[t,T]$. Using \eqref{leibniz2} formally (this will later be justified by the regularity properties of the involved 
functions), we then obtain for every $t\in [0,T]$ the
representation formula
\gianni{%
\begin{align}
\label{meanpq}
\mbox{mean}\,(\pal(t)+\tau\qal(t))\,&= \,\mbox{mean}\,(g_1^\alpha)\,+\,\frac 1{|\Omega|}
\int_t^T\!\!\iO(g_2^\alpha-\psi_1^\alpha\,\qal-\psi_2^\alpha\,\qal),
\end{align}
}%
where the \lhs\ equals $\,\mbox{mean}\,(\pal(t))\,$ for almost every $\,t\in (0,T)$ \gianni{by~\eqref{meanqal}}. 
In view of this identity,
we cannot expect the bound \eqref{lgq} to hold also in this case: indeed, due to the presence of the 
term  $\,-\int_t^T\!\!\iO\psi_2^\alpha\,\qal$ on the \rhs\ of \eqref{meanpq}, we cannot hope
to be able to control the mean value of $\,\pal$ independently of $\alpha>0$.

Nevertheless, a proper solution to \eqref{adj1}--\eqref{adj3} exists also in this case.
To this end, we eliminate $\,\mbox{mean}\,(\pal)\,$ from the problem, following a strategy
 introduced in \cite{CFGS1} and~\cite{CGSAMO}. We put
 \begin{equation}\label{defH0}
 H_0:=\{v\in H: \mbox{\,mean}(v)=|\Omega|^{-1}(v,{\mathbf{1}})=0\}.
 \end{equation}
 Then $\,H=H_0\oplus \mbox{span}\{{\mathbf{1}}\}$, and we have (cf.~\eqref{defVrpos}) that
  $\,V_0^r=V_A^r\cap H_0\,$ for $\lambda_1=0$. Moreover, the
 linear operator $\,A_0^{2r}=A^{2r}_{|V_0^{2r}}\,$ is a topological isomorphism from 
 $V_0^{2r}$ onto~$H_0$, where we have the representation formulas
\begin{align}
 \label{arnull}
 &A_0^{2r}v\,=\,A^{2r}v\,=\,\sum_{j=2}^\infty \lambda_j^{2r}(v,e_j)e_j \quad\forall \,v\in V_0^{2r},\\
 \label{a-rnull}
 &A_0^{-2r}v\,:=\,(A_0^{2r})^{-1}v\,=\, \sum_{j=2}^\infty \lambda_j^{-2r}(v,e_j)e_j \quad\forall \,v\in H_0.
\end{align} 
Moreover, with
\Beq\label{defh10}
H^{1,0}(\Omega):=\Huno \cap H_0,
\Eeq
we have (cf.~\cite[Sect.~5]{CGS19}) that $\,(H^{1,0}(\Omega),H_0, (H^{1,0}(\Omega))^*)\,$ is a 
Hilbert triple with dense, continuous,  and compact embeddings.
 
Now observe that\, $A^r(\mbox{mean}\,(\pal(t))\mathbf{1})=\mbox{mean}\,(\pal(t))\, A^r{\mathbf{1}}=0$, 
and thus \gianni{\eqref{adj1} becomes}
\begin{align}\label{ll5}
&(A^r(\pal(t)-\mbox{mean}(\pal(t))\mathbf{1}),A^r v)\,=\,(\qal(t),v) \non\\
&\quad\pier{\mbox{for a.e. \,$t\in (0,T)$\, and every $v\in V_A^r$}.}
\end{align}   
Since $\,\pal(t)-\mbox{mean}(\pal(t))\mathbf{1}\in H_0\,$ for almost every $t\in (0,T)$, this is equivalent to
\begin{equation}
\label{ll6}
A_0^{2r}(\pal-\mbox{mean}(\pal)\mathbf{1})=\qal \quad\mbox{and}\quad \pal-\mbox{mean}(\pal)\mathbf{1}\,=\,A_0^{-2r}\qal .
\end{equation}

At this point, we are able to state the existence result for the system \eqref{adj1}--\eqref{adj2}
in the case $\,\lambda_1=0$ by adapting the results established in \cite[Sect.~5]{CGS19} to the
present situation. We then can infer that there
exists a unique solution $(\pal,\qal)$ such that
\begin{align}
\label{hugo1}
&A_0^{-2r}\qal\in L^\infty(0,T;V_0^{2r}),\\[1mm]
\label{hugo2}
&\qal\in L^\infty(0,T;H_0)\cap L^2(0,T;H^{1,0}(\Omega)),\\[1mm]
\label{hugo3}
&(A_0^{-2r}+\tau I)\qal\in H^1(0,T;(H^{1,0}(\Omega))^*),
\end{align}
as well as 
\begin{align}
&\mbox{mean}(\pal+\tau\qal) \,\mbox{ satisfies }\,\eqref{meanpq}\,\mbox{ for every }\,t\in [0,T],\\[2mm] 
\label{hugo5}
&\pal-\mbox{mean}(\pal)\mathbf{1}=A_0^{-2r}\qal,\\[1mm]
\label{hugo6}
&\left\langle -\dt\bigl(A_0^{-2r}+\tau I\bigr)\qal(t),v\right\rangle_{H^{1,0}(\Omega)}
\,+\,\iO\nabla\qal(t)\cdot\nabla v \,+\,\bigl((\psi_1^\alpha(t)+\gianni{\psi_2^\alpha}(t))\,\qal(t),v\bigr)
\non\\[1mm]
&\quad =\,(g_2^\alpha(t),v)\quad\mbox{for a.e. \,$t\in (0,T)$\, and every }\,v\in H^{1,0}(\Omega),
\\[2mm]
\label{hugo7}
&\left\langle (A_0^{-2r}+\tau I)\qal(T),v \right\rangle_{H^{1,0}(\Omega)}
\,=\,\bigl(g_1^\alpha-\mbox{mean}(g_1^\alpha)\mathbf{1},v)\quad\mbox{for all }\,v\in H^{1,0}
(\Omega).
\end{align}
Notice that, by \eqref{hugo3}, we have $\,(A_0^{-2r}+\tau I)\qal \in C^0([0,T];(H^{1,0}(\Omega))^*)$,
which gives the endpoint condition \eqref{hugo7} a proper meaning: indeed, \eqref{hugo7}
means that $\,(A_0^{-2r}+\tau I)\qal(T)=g_1^\alpha-\mbox{mean}(g_1^\alpha)\mathbf{1}\,$ in
$\,(H^{1,0}(\Omega))^*$, where the \rhs\ belongs to~$H_0$. Now observe that the operator
\Beq
\label{hugo8}
(A_0^{-2r}+\tau I)v\,=\,\sum_{j=2}^\infty (\lambda_j^{-2r}+\tau)(v,e_j)e_j \quad\forall \,v\in H_0
\Eeq  
is a topological isomorphism from $H_0$ into itself with the inverse
\Beq
\label{hugo9}
(A_0^{-2r}+\tau I)^{-1}v\,=\,\sum_{j=2}^\infty (\lambda_j^{-2r}+\tau)^{-1}(v,e_j)e_j 
\quad\forall \,v\in H_0.
\Eeq
Hence, also $\,\qal(T)=(A_0^{-2r}+\tau I)^{-1}(g_1^\alpha-\mbox{mean}(g_1^\alpha)\mathbf{1})\,$
has a proper meaning as an element of~$H_0$.
Next, we consider the mapping
\begin{align} 
 \label{hugo10}
 &(A_0^{-2r}+\tau I)^{1/2}v\,=\,\sum_{j=2}^\infty (\lambda_j^{-2r}+\tau)^{1/2}(v,e_j)e_j\quad\forall\,v\in H_0.
 \end{align}
It is readily seen that the estimates \eqref{lg5}--\eqref{lg7} have the analogues
\begin{align}
\label{hugo11}
&\|(A_0^{-2r}+\tau I)^{1/2}v\|^2\,\ge\,\tau\,\|v\|^2\quad\forall\,v\in H_0,\\[1mm]
\label{hugo12}
&\|(A_0^{-2r}+\tau I)^{1/2}v\|^2\,\le\,\gianni{(\lambda_2^{-2r}+\tau)}\,\|v\|^2\quad\forall\,v\in H_0,\\[1mm] 
\label{hugo13}
&\|(A_0^{-2r}+\tau I)^{-1}v\|^2\,\le\,\tau^{-2}\,\|v\|^2\quad\forall\,v\in H_0.
\end{align}

Now observe that for a.e. $t\in(0,T)$ it holds that
\begin{equation}
\label{hugo14}
-\langle (A_0^{-2r}+\tau I)\,\qal(t),\qal(t)\rangle_{H^{1,0}(\Omega)}\,=
\,-\,\frac 12\,\frac d{dt}\left\|(A^{-2r}_0+\tau I)^{1/2}
\qal(t)\right\|^2 .
\end{equation}
At this point, we insert $\,v=\qal(t)\in H^{1,0}(\Omega)$ in \eqref{hugo6} and integrate over $[t,T]$, where $t\in [0,T)$, to recover the
identity \eqref{lg3}, only that in the first term the expression $A^{-2r}$ and 
the dual pairing between $\Huno^*$ and $\Huno$ are replaced by $A_0^{-2r}$ and the
dual pairing between $(H^{1,0}(\Omega))^*$ and $H^{1,0}(\Omega)$, respectively.
Again, the third summand on the \lhs\ is nonnegative, and the \rhs\ is bounded by the expression
\eqref{Paul}. Moreover,  the first summand on the \lhs, 
which we denote by $I_1^\alpha(t)$, can by \eqref{hugo12} and \eqref{hugo13} be estimated as follows:
\begin{align}
\label{hugo15}
I_1^\alpha(t)&\,=\,\frac 12 \left\|(A_0^{-2r}+\tau I)^{1/2}\qal(t)\right\|^2 \,-\,
\frac 12 \left\|(A_0^{-2r}+\tau I)^{1/2}\qal(T)\right\|^2 \,\ge\,
\frac \tau 2\,\|\qal(t)\|^2-\frac 12\,C_6.
\end{align}
At this point, we can again employ Gronwall's lemma  to conclude the
estimate 
\begin{equation}
\label{llq}
\|\pal-\mbox{mean}(\pal)\mathbf{1}\|_{L^\infty(0,T;V_0^{2r})}\,+\,
\|\qal\|_{L^\infty(0,T;H_0)\cap L^2(0,T;H^{1,0}(\Omega))}\,
\le\,C_7 \quad\forall\,\juerg{\alpha\in (0,1]},
\end{equation}
which is the sought analogue of \eqref{lgq}. 

\vspace{3mm}
In the following, we complement \eqref{lgq} and \eqref{llq} by further estimates.
We treat the two cases $\lambda_1>0$ and $\lambda_1=0$ simultaneously, where it is understood that 
the spaces $V_0^r$ and the operators $A_0^r$ are defined as in \eqref{defVrpos} and \eqref{defAz}, respectively.
We now introduce the space   
\begin{equation}
\label{defZ}
{\cal Z}:=\left\{
\begin{array}{ll}
\{v\in H^1(0,T;H)\cap L^2(0,T;\Huno):v(0)=0\}&\mbox{if $\,\lambda_1>0$}\\[1mm]
\{v\in H^1(0,T;H_0)\cap L^2(0,T;H^1(\Omega)):v(0)=0\}&\mbox{if $\,\lambda_1=0$}
\end{array}
\right. ,
\end{equation}
which is a Hilbert space when endowed with its natural inner product and norm. 
Moreover, setting 
\begin{equation}
\label{defG}
G=H \quad\mbox{for }\,\lambda_1>0\quad\mbox{and}\quad G=H_0\quad\mbox{for }\, \lambda_1=0,
\end{equation}
we see that the embedding $\,{\cal Z}\subset C^0([0,T];G)\,$ is continuous. 
\gianni{\pier{Furthermore,} we also have the dense and continuous embeddings ${\cal Z}\subset L^2(0,T;G)\subset{\cal Z}^*$}, 
where it is understood that 
\begin{equation}
\label{dualZ}
\langle v,z\rangle_{\cal Z}=\int_0^T(v(t),z(t))\,dt \quad\mbox{for all $\,z\in {\cal Z}\,$ and \,$v\in 
L^2(0,T;G)$.}
\end{equation}

In order to avoid to have to distinguish between the two cases, we employ in the following
the same notation $\,\langle\cdot,\cdot\rangle\,$ for the dual pairings 
$\,\langle\cdot,\cdot\rangle_{\Huno}\,$  and $\,\langle\cdot,\cdot\rangle_{H^{1,0}(\Omega)}$,
where the former corresponds to the case $\lambda_1>0$ and the latter to the case $\lambda_1=0$. 

At this point, \gianni{by recalling \accorpa{regp}{regpq} for $\lambda_1>0$ and \accorpa{hugo1}{hugo3} for $\lambda_1=0$},
we may employ the integration-by-parts formula \eqref{leibniz2} \gianni{with $z=A_0^{-2r}\qal+\tau\qal$}
to conclude that for every $v\in {\cal Z}$ it holds that
\begin{align}
\label{absch1}
\langle -\dt(\pal+\tau\qal\pier{)},v\rangle_{{\cal Z}}\,&=\,-\int_0^T\langle \dt(A_0^{-2r}\qal(t)+\tau\qal(t),v(t)\rangle\,dt
\non\\
&=\,\int_0^T(\dt v(t),(A_0^{-2r}+\tau I)\qal(t))\,dt\,-\,(g_1^\alpha,v(T))\non\\[1mm]
&\le\,\|\dt v\|_{L^2(0,T;H)}\,\|(A_0^{-2r}+\tau I)\qal\|_{L^2(0,T;H)}\,+\,\|g_1^\alpha\|_H\,\|v(T)\|_H\non\\[1mm]
&\le\,C_7\,\|v\|_{\cal Z},
\end{align}
which implies that
\begin{equation}
\label{absch2}
\|\dt(\pal+\tau\qal)\|_{{\cal Z}^*}\,\le\,C_7 \quad\forall\,\juerg{\alpha\in(0,1]}.
\end{equation}
Now observe that for any $v\in {\cal Z}$ it holds that
\begin{align*}
&\int_0^T(\nabla\qal(t),\nabla v(t))\,dt\,+\int_0^T(\psi_1^\alpha(t)\qal(t),v(t))\,dt\,
-\int_0^T (g_2^\alpha(t),v(t))\,dt
\\[1.5mm]
&\le\,\|\qal\|_{L^2(0,T;\Huno)}\,\|v\|_{\cal Z}\,+\,C_8\,\|\qal\|_{L^2(0,T;H)}\,\|v\|_{L^2(0,T;H)}\,
+\,C_9\,\|v\|_{L^2(0,T;H)}\\[1.5mm]
&\le\,C_{10}\,\|v\|_{\cal Z},
\end{align*}
and it follows from comparison in \eqref{adj2} that, with $\,\Lambda^\alpha:=\psi_2^\alpha \qal=
\phi(\alpha)h''(\byal)\qal$,
\begin{equation}
\label{lamal}
\|\Lambda^\alpha\|_{{\cal Z}^*}\,\le\,C_{11} \quad\forall\,\juerg{\alpha\in(0,1]}.
\end{equation}

At this point, we choose any sequence $\{\alpha_n\}$ such that $\alpha_n\searrow0$. We infer from Theorem 3.3 
and Theorem 4.5 that, at least for a subsequence which is again indexed by~$\,n$, 
\begin{align}
\label{ustark}
&{\bu}^{\alpha_n}\to\bu\quad\mbox{strongly in }\,L^2(Q),\\[0.5mm]
\label{yweak}
&{\by}^{\alpha_n}\to\by \quad\mbox{weakly-star in }\,W^{1,\infty}(0,T;H)\cap H^1(0,T;\Huno).
\end{align}
By virtue of \cite[Sect.~8,~Cor.~4]{Simon}, we may also assume that
\begin{align}
\label{stark1}
&{\by}^{\alpha_n}\to\by \quad\mbox{strongly in }\,C^0([0,T];L^p(\Omega)) \quad\forall\,p\in [1,6),
\end{align}
which entails, in particular, that
\begin{align}
\label{stark2}
&f_2''({\by}^{\alpha_n})\to f_2''(\by) \quad\mbox{strongly in }\,C^0([0,T];L^p(\Omega)) \quad\forall
\,p\in [1,6),\\[0.5mm]
\label{stark3}
&g_1^{\alpha_n}\to \beta_1(\by(T)-y_\Omega)\quad\mbox{strongly in }\,H,\\[0.5mm]
\label{stark4}
&g_2^{\alpha_n}\to \beta_2(\by-y_Q)\quad\mbox{strongly in }\,L^2(Q).
\end{align} 
Moreover, by virtue of the estimates \eqref{lgq}, \eqref{llq}, \eqref{absch2}, and \eqref{lamal},
there are limits $\,\zeta,\bar q,\Lambda$ such that, at least for another subsequence which is again
indexed by $n$,
\begin{align}
\label{limzeta}
\dt(A_0^{-2r}+\tau I)q^{\alpha_n}&\to\zeta\quad\mbox{weakly in }\,{\cal Z}^*,\\[0.5mm]
\label{limq}
q^{\alpha_n}&\to \bar q\quad\mbox{weakly-star in }\,L^\infty(0,T;G)\cap L^2(0,T;\Huno),\\[0.5mm]
\label{limA0q}
A_0^{-2r}q^{\alpha_n}&\to A_0^{-2r}\bar q\quad\mbox{weakly-star in }\,L^\infty(0,T;V_0^{2r}),\\[0.5mm] 
\label{limlam}
\Lambda^{\alpha_n}&\to\Lambda\quad\mbox{weakly in }\,{\cal Z}^*.
\end{align}
The limit $\zeta\in {\cal Z}^*$ is readily identified. Indeed, by formula \eqref{leibniz2} we have,
for every $v\in{\cal Z}$,
\begin{align}
\label{findzeta}
&\lim_{n\to\infty}\int_0^T\langle\dt(A_0^{-2r}+\tau I)q^{\alpha_n}(t),v(t)\rangle\,dt\non\\
&=\,\lim_{n\to\infty}\Bigl[-\int_0^T\bigl(\dt v(t),(A_0^{-2r}+\tau I)q^{\alpha_n}(t)\bigr)\,dt
+\,(g_1^{\alpha_n},v(T))\Bigr]\non\\
&=\,-\int_0^T\bigl(\dt v(t),(A_0^{-2r}+\tau I)\bar q(t)\bigr)\,dt\,+\,(\beta_1(\by(T)-y_\Omega)
,v(T))\,=:\,\langle\zeta,v\rangle_{{\cal Z}}.
\end{align} 
Moreover, by combining the strong convergence \eqref{stark2} with \eqref{limq}, it is easily checked that
\begin{equation}
\label{yeah}
f_2''({\by}^{\alpha_n})\,q^{\alpha_n}\to f_2''(\by)\,\bar q\quad\mbox{weakly in }\,L^2(Q).
\end{equation} 
At this point, we recall that 
$$
\langle -\dt(\pal+\tau\qal)(t),v(t)\rangle\,=\,\langle -\dt(A_0^{-2r}+\tau I)\qal(t),v(t)\rangle
\quad\mbox{for a.e.$\,t\in (0,T)\,$ and \,$v\in{\cal Z}$.}
$$ 
We now write the adjoint system \eqref{adj1}--\eqref{adj2} for $\alpha=\alpha_n$, insert $v=v(t)$
for an arbitrary $v\in {\cal Z}$, integrate the resulting identity with respect to time over $[0,T]$, and pass to 
the limit as $n\to\infty$. It then results the following equation:
\begin{align}
\label{limeq}
\langle \Lambda,v\rangle_{\cal Z}\,&=\,-\int_0^T\bigl(\dt v(t),(A_0^{-2r}+\tau I)\bar q(t)\bigr)\,dt
\,+\beta_1(\by(T)-y_\Omega,v(T))\,\non\\
&\quad\,\,\pier{{}-\int_0^T\!\!\iO\nabla\bar q\cdot\nabla v}\,+\int_0^T\!\!\iO\bigl(\beta_2(\by-y_Q)-f_2''(\by)\bar q\bigr)\,v \quad\forall\,v\in{\cal Z}.
\end{align}

Finally, we need to identify the variational inequality relating the optimal control to the adjoint variables. 
In this regard, we can infer, with the same argument as in the proof of \cite[Thm.~5.9]{CGS19}, that
the optimal control $\bu^{\alpha_n}$ satisfies the variational inequality
\begin{align}\label{vugtil}    
\int_0^T\!\!\iO(q^{\alpha_n}+\beta_3\bar u^{\alpha_n} + \bar u^{\alpha_n}-\bu)
(v-\bar u^{\alpha_n})\,\ge\,0 \quad\forall\,v\in\uad.
\end{align}
Taking the limit as $n\to\infty$ in \eqref{vugtil}, and using \eqref{ustark} and \gianni{\eqref{limq}}, we
arrive at the necessary optimality condition
\begin{align}
\label{vugo}
\int_0^T\!\!\int_\Omega (\bar q\,+\,\beta_3\,\bu)(v-\bu)\,\ge\,0\quad\forall\,v\in\uad. 
\end{align}

From the above considerations, we can conclude the following first-order necessary
\gianni{optimality} conditions for the optimal control problem~\CP0:

\Bthm
Suppose that the conditions \pier{{\bf (A1)}--{\bf {(A6)}},  {\bf (A8)}, {\bf {(A9)}}} are satisfied, and let $\,\bu\in\uad\,$ be an optimal control
for {\rm (${\mathcal{CP}}_0$)} with associated state $\,(\bm,\by)\,$ where $\,\by=\S_0(\bu)$. 
Then there exist $(\bar q,\Lambda)$ such that the 
following statements hold true:

\noindent
{\rm (i)} We have the regularity properties
\begin{align}
\bar q\in L^\infty(0,T;G)\cap L^2(0,T;\Huno), \quad \Lambda\in {\cal Z}^*.
\end{align}

\noindent
{\rm (ii)} \,\,The adjoint equation \eqref{limeq} is fulfilled.

\noindent
{\rm (iii)} \,The necessary optimality condition \eqref{vugo} is satisfied.
\Ethm

\Brem
From \eqref{vugo} we infer that\pier{, in the case $\beta_3>0$, $\,\bu\,$ is} nothing but the
$L^2(Q)$-orthogonal projection of \,$-\beta_3^{-1}\,q\,$ onto $\uad$.
\Erem
\Brem
Unfortunately, we are unable to derive any complementarity slackness 
conditions for the Lagrange multiplier $\Lambda$. Indeed, while it is easily seen that
$$
\liminf_{n\to\infty}\int_0^T\!\!\iO\Lambda^{\alpha_n}\,q^{\alpha_n}\,=\,\liminf_{n\to\infty}\int_0^T\!\!\iO\frac 
{2 \phi(\alpha_n)}{1-(\bar y^{\alpha_n})^2}\,|q^{\alpha_n}|^2\,\ge\,0
\quad \forall\,n\in\enne,
$$
the convergence properties \eqref{yweak} and \eqref{limq} do not suffice to conclude that
$\,\langle \Lambda,\bar q\rangle_{\cal Z}\,\ge\,0$.   
\Erem


\section*{Acknowledgments}
\pier{This research was supported by the Italian Ministry of Education, 
University and Research~(MIUR): Dipartimenti di Eccellenza Program (2018--2022) 
-- Dept.~of Mathematics ``F.~Casorati'', University of Pavia. 
In addition, PC and CG gratefully acknowledge some other 
financial support from the MIUR-PRIN Grant 2015PA5MP7 ``Calculus of Variations'',}
the GNAMPA (Gruppo Nazionale per l'Analisi Matematica, 
la Probabilit\`a e le loro Applicazioni) of INdAM (Isti\-tuto 
Nazionale di Alta Matematica) and the IMATI -- C.N.R. Pavia.


\vspace{3truemm}

{\small

\Begin{thebibliography}{10}

\bibitem{AM}
M. Ainsworth, Z. Mao, Analysis and approximation of a fractional
Cahn--Hilliard equation. SIAM J. Numer. Anal. {\bf 55} (2017),
1689-1718.

\bibitem{AkSeSchi}
G. Akagi, G. Schimperna, A. Segatti,
Fractional Cahn--Hilliard, Allen--Cahn and porous medium equations.
J. Differ. Equations {\bf 261} (2016), 2935-2985.

\juerg{
\bibitem{AKW}
H. Antil, R. Khatri, M. Warma, External optimal control of nonlocal PDEs.
Preprint arXIv:1811.04515 [math.OC] (2018), pp. 1-30.
}

\bibitem{HO1}
H. Antil, E. Ot\'arola,
A FEM for an optimal control problem of fractional
powers  of  elliptic  operators. SIAM J. Control Optim. {\bf 53} (2015), 3432-3456.

\bibitem{HO2}
H. Antil, E. Ot\'arola, An a posteriori error analysis 
for an optimal control problem involving the fractional Laplacian.
IMA J. Numer. Anal. {\bf 38} (2018), 198-226.

\bibitem{HAS}
H. Antil, E. Ot\'arola, A.J. Salgado, A space-time fractional optimal control problem:
analysis and discretization.
 SIAM J. Control Optim. {\bf  54} (2016), 1295-1328. 

\juerg{
\bibitem{APR}
H. Antil, J. Pfefferer, S. Rogers, Fractional operators with inhomogeneous boundary conditions:
analysis, control and discretization. Preprint arXiv:1703.05256 [math.NA] (2017), pp. 1-38.%
}

\juerg{%
\bibitem{AW1}
H. Antil, M. Warma, Optimal control of fractional semilinear PDEs. Preprint arXiv:1712.04336
[math.OC] (2017), pp. 1-27.%
}

\juerg{
\bibitem{AW2}
H. Antil, M. Warma, Optimal control of the coefficient for the regional fractional $p$--Laplace equation: approximation and convergence. AIMS Math. Control Related Fields (2018),
https://doi.org/10.3934/mcrf.2019001. 
}

\bibitem{Barbu}
V. Barbu,
``Nonlinear Differential Equations of Monotone Type in Banach Spaces''.
Springer, London, New York, 2010.

\bibitem{BDM1}
T. Biswas, S. Dharmatti,  M.T. Mohan, Pontryagin's maximum principle for optimal
control of the nonlocal Cahn--Hilliard--Navier--Stokes systems in two dimensions.
Preprint arXiv:1802.08413 [math.OC] (2018), pp.~1-36.
 
\bibitem{BDM2}
T. Biswas, S. Dharmatti, M.T. Mohan, Maximum principle and data assimilation problem
for the optimal control problems governed by 2D nonlocal Cahn--Hilliard--Navier--Stokes equations.
Preprint arXiv:1803.11337 [math.AP] (2018), pp.~1-26.

\bibitem{CFGS1}
P. Colli, M.H. Farshbaf-Shaker, G. Gilardi, J. Sprekels,
Optimal boundary control of a viscous Cahn--Hilliard system with dynamic boundary condition 
and double obstacle potentials. SIAM J.
Control Optim. {\bf 53} (2015), 2696-2721.

\bibitem{CFGS2}
P. Colli, M.H. Farshbaf-Shaker, G. Gilardi, J. Sprekels,
Second-order analysis of a boundary
control problem for the viscous Cahn--Hilliard equation with dynamic boundary conditions.
Ann. Acad. Rom. Sci. Ser. Math. Appl. {\bf 7}
(2015), 41-66.

\bibitem{CFGS3}
P. Colli, M.H. Farshbaf-Shaker, J. Sprekels,
A deep quench approach to the optimal control of
an Allen--Cahn equation with dynamic boundary conditions and double obstacles.
Appl. Math. Optim. {\bf 71}
(2015), 1-24.

\bibitem{CGRS}
P. Colli, G. Gilardi, E. Rocca, J. Sprekels,
Optimal distributed control of a diffuse interface model
of tumor growth.
Nonlinearity {\bf 30}
(2017), 2518-2546.

\bibitem{CGSAdvan}
P. Colli, G. Gilardi, J. Sprekels,
A boundary control problem for the pure Cahn--Hilliard equation
with dynamic boundary conditions.
Adv. Nonlinear Anal. {\bf 4}
(2015), 311-325.

\bibitem{CGSAIMS}
P. Colli, G. Gilardi, J. Sprekels, 
Distributed optimal control of a nonstandard nonlocal phase field
system. AIMS Mathematics {\bf 1}
(2016), 246-281.

\bibitem{CGSAMO}
P. Colli, G. Gilardi, J. Sprekels,
A boundary control problem for the viscous Cahn--Hilliard equation
with dynamic boundary conditions.
Appl. Math. Optim. {\bf 72}
(2016), 195-225.

\bibitem{CGSEECT}
P. Colli, G. Gilardi, J. Sprekels,
Distributed optimal control of a nonstandard nonlocal phase field
system with double obstacle potential.
Evol. Equ. Control Theory {\bf 6}
(2017), 35-58.

\bibitem{CGSAnnali}
P. Colli, G. Gilardi, J. Sprekels, On a Cahn--Hilliard system with convection
and dynamic boundary conditions. Ann. Mat. Pura Appl.~(4) {\bf 197}
(2018), 1445-1475.

\bibitem{CGSSIAM18}
P. Colli, G. Gilardi, J. Sprekels,
Optimal velocity control of a viscous Cahn--Hilliard system with convection and dynamic boundary conditions.
SIAM J. Control Optim. {\bf 56} (2018), 1665-1691.

\bibitem{CGSConvex}
P. Colli, G. Gilardi, J. Sprekels,
Optimal velocity control of a convective Cahn--Hilliard system with double
obstacles and dynamic boundary conditions: a `deep quench' approach.
\pier{To appear in J. Convex Anal. {\bf 26} (2019) (see also Preprint arXiv: 1709.03892
[math.AP] (2017), pp.~1-30).}

\bibitem{CGS18}
P. Colli, G. Gilardi, J. Sprekels, 
Well-posedness and regularity for a generalized fractional Cahn--Hilliard system. 
Preprint arXiv:1804.11290 [math.AP] (2018), \pier{pp.~1-36}.

\bibitem{CGS19}
P. Colli, G. Gilardi, J. Sprekels, 
Optimal distributed control of a generalized fractional Cahn--Hilliard system.
\juerg{Appl. Math. Optim. (2018), https://doi.org/10.1007/s00245-018-9540-7.}

\bibitem{CS}
P. Colli, J. Sprekels, Optimal boundary control of a nonstandard Cahn--Hilliard system 
with dynamic boundary condition and double obstacle inclusions. In:
``Solvability, Regularity, \juerg{and}
Optimal Control of Boundary Value Problems for PDEs''  (P. Colli, A. Favini, E. Rocca, 
G. Schimperna, J. Sprekels, eds.), Springer INdAM Series Vol. {\bf 22}, pp. 151-182,
Springer 2017.

\bibitem{Duan}
N. Duan, X. Zhao,
Optimal control for the multi-dimensional viscous Cahn--Hilliard equation. 
Electron. J. Differ. Equations 2015, Paper No. 165, 13 pp.


\bibitem{FGG}
S. Frigeri, C.G. Gal, M. Grasselli, 
On nonlocal Cahn--Hilliard--Navier--Stokes systems in two dimensions. 
J. Nonlinear Sci. {\bf 26}
(2016), 847-893.

\bibitem{FGGS}
S. Frigeri, C.G. Gal, M. Grasselli, J. Sprekels, 
\juerg{Two-dimensional nonlocal Cahn--Hilliard--Navier--Stokes
systems with variable viscosity, degenerate mobility and singular potential. To appear in
Nonlinearity (see also WIAS Preprint Series No. 2309, Berlin 2016, \pier{pp.~1-56}).}

\bibitem{FGS}
S. Frigeri, M. Grasselli, J. Sprekels, Optimal distributed control of two-dimensional 
nonlocal Cahn--Hilliard--Navier--Stokes systems with degenerate mobility and singular potential.
\juerg{Appl. Math. Optim. (2018), https://doi.org/10.1007/s00245-018-9524-7.}

\bibitem{FRS}
S. Frigeri, E. Rocca, J. Sprekels,
Optimal distributed control of a nonlocal Cahn--Hilliard/Navier--Stokes
system in two dimensions. SIAM J. Control Optim.
{\bf 54} (2016), 221-250.

\bibitem{Fukao}
T. Fukao, N. Yamazaki, A boundary control problem for the equation and dynamic
boundary condition of Cahn--Hilliard type. In: ``Solvability, Regularity, \juerg{and}
Optimal Control of Boundary Value Problems for PDEs''  (P. Colli, A. Favini, E. Rocca, 
G. Schimperna, J. Sprekels, eds.), Springer INdAM Series Vol. {\bf 22}, pp. 255-280, Springer 2017.

\bibitem{GalDCDS}
C.G. Gal, 
On the strong-to-strong interaction case 
for doubly nonlocal Cahn--Hilliard equations.
Discrete Contin. Dyn. Syst.  {\bf 37} (2017),  131-167.
\bibitem{GalEJAM}
C.G. Gal, 
Non-local Cahn--Hilliard equations with fractional dynamic boundary conditions. European J. Appl. Math. {\bf 28} (2017), 736-788. 
\bibitem{GalAIHP}
C.G. Gal, 
Doubly nonlocal Cahn--Hilliard equations. 
Ann. Inst. H. Poincar\'e Anal. Non Lin\'eaire {\bf 35} (2018), 357-392.

\bibitem{GV}
C. Geldhauser, E. Valdinoci,
Optimizing the fractional power in a model with stochastic PDE constraints.
Preprint 	arXiv:1703.09329v1 [math.AP] (2017), pp. 1-18.

\bibitem{HHCK}
M. Hinterm\"uller, M. Hinze, C. Kahle, T. Keil,
A goal-oriented dual-weighted adaptive finite element approach for the optimal control 
of a nonsmooth Cahn--Hilliard--Navier--Stokes system.
Optimization and Engineering (2018), https://doi.org/10.1007/s11081-018-9393-6.

\bibitem{HKW}
M.  Hinterm\"uller,  T.  Keil,  D.  Wegner,
Optimal  control  of  a  semidiscrete  Cahn--Hilliard--Navier--Stokes
system with non-matched fluid densities. SIAM J. Control Optim. {\bf 55} (2018),
1954-1989.

\bibitem{HW1}
M. Hinterm\"uller, D. Wegner,
Distributed optimal control of the Cahn--Hilliard system including the case
of a double obstacle homogeneous free energy density. SIAM J. Control Optim.
{\bf 50} (2012), 388-418.

\bibitem{HW2}
M. Hinterm\"uller, D. Wegner,
Optimal control of a semidiscrete Cahn--Hilliard--Navier--Stokes system.
SIAM J. Control Optim. {\bf 52} (2014), 747-772.

\bibitem{HW3}
 M. Hinterm\"uller, D. Wegner,
Distributed and boundary control problems for the semidiscrete Cahn--Hilliard/Navier--Stokes 
system with nonsmooth Ginzburg--Landau energies. In: ``Topological Optimization
and Optimal Transport'', Radon Series on Computational and Applied Mathematics vol.
{\bf 17} (2017), pp. 40-63.

\bibitem{Medjo}
T. Tachim Medjo, 
Optimal control of a Cahn--Hilliard--Navier--Stokes model 
with state constraints.
J. Convex Anal.
{\bf 22} (2015), 1135-1172.

\bibitem{RS}
E. Rocca, J. Sprekels,
Optimal distributed control of a nonlocal convective Cahn--Hilliard equation by
the velocity in three dimensions. SIAM J. Control Optim. {\bf 53} (2015), 1654-1680.

\bibitem{Simon}
J. Simon,
Compact sets in the space $L^p(0,T; B)$.
 Ann. Mat. Pura Appl.~(4) 
{\bf 146}, (1987), 65-96.

\bibitem{SV}
J. Sprekels, E. Valdinoci, A new class of identification problems: optimizing the
fractional order in a nonlocal evolution equation. SIAM J. Control Optim. {\bf 55}
(2017), 70-93.

\bibitem{WN}
Q. F. Wang, S.-i. Nakagiri, 
Weak solutions of Cahn--Hilliard equations having forcing terms and
optimal control problems. Mathematical models in functional equations 
(Kyoto, 1999), Surikaisekikenkyusho Kokyuroku No. 1128 (2000), 172-180.

\bibitem{ZL1}
 X.P. Zhao, C.C. Liu,
Optimal control of the convective Cahn--Hilliard equation. Appl. Anal.
{\bf 92} (2013), 1028-1045.

\bibitem{ZL2}
X.P. Zhao, C.C. Liu,
Optimal control of the convective Cahn--Hilliard equation in 2D case. Appl.
Math. Optim. {\bf 70}
(2014), 61-82.

\bibitem{Z}
J. Zheng, Time optimal controls of the Cahn--Hilliard equation with 
internal control. Optimal Control Appl. Methods {\bf 36}
(2015), 566-582.

\bibitem{ZW}
J. Zheng, Y. Wang, Optimal control problem for Cahn--Hilliard equations with 
state constraint. J. Dyn. Control Syst. {\bf 21} (2015), 257-272.

\End{thebibliography}}

\End{document}
